\documentclass[11pt,a4paper]{article}
\usepackage[T1]{fontenc}
\usepackage[utf8]{inputenc}
\usepackage{authblk}
\newcommand{\mycomment}[1]{}
\usepackage{color}

\usepackage{graphicx}
\graphicspath{{./figures/}}
\everymath{\displaystyle}
\usepackage{booktabs}
\usepackage{enumitem}
\usepackage{amsmath}
\usepackage{amssymb}
\usepackage{xspace}
\usepackage{latexsym}
\usepackage{epsfig}
\usepackage{verbatim}
\usepackage{alltt}
\usepackage{url}
\usepackage{longtable}
\usepackage{graphicx}
\usepackage{subfigure}
\usepackage{epstopdf}
\usepackage{float}
\usepackage[T1]{fontenc}
\usepackage[utf8]{inputenc}
\usepackage[font=small,labelfont=bf]{caption}

\usepackage{slashbox}
\newtheorem{thm}{Theorem}

\newcommand{\bbb}{\beta}

\newcommand{\ba}{\begin{array}} 
\newcommand{\ea}{\end{array}}
\newcommand{\be}{\begin{equation}} 
\newcommand{\ee}{\end{equation}}
\newcommand{\bea}{\begin{eqnarray}}
\newcommand{\eea}{\end{eqnarray}}
\newcommand{\bean}{\begin{eqnarray*}}
\newcommand{\eean}{\end{eqnarray*}}

\allowdisplaybreaks

\newcommand{\Keywords}[1]{\par\noindent
{\small{\em Keywords\/}: #1}}

\title{Numerical Implementation of a  Cohesive Zone Model in History-Dependent Materials}

\author{L. Hakim\thanks{layal.hakim@brunel.ac.uk}\ }
\author{S.E. Mikhailov\thanks{sergey.mikhailov@brunel.ac.uk (Correponding author)}}
\affil{Department of Mathematical Sciences, Brunel University, London, UK}

\begin{document}
\maketitle
\begin{abstract}
A non-linear history-dependent cohesive zone model of crack propagation in linear elastic and visco-elastic materials is presented. The viscoelasticity is described by a linear Volterra integral operator in time. The normal stress on the cohesive zone satisfies the history dependent yield condition, given by a non-linear Abel-type integral operator.
The crack starts propagating, breaking the cohesive zone, when the crack tip opening reaches a prescribed critical value. A numerical algorithm for computing the evolution of the crack and cohesive zone in time is discussed along with some numerical results.

\end{abstract}
\section{Introduction}
\begin{figure}[H]
\begin{center}
\includegraphics[scale=0.43]{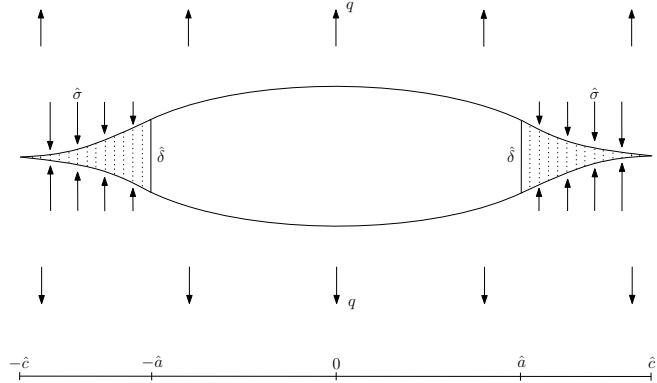}
\caption{Cohesive zone}\label{cohesivepic}
\end{center}
\end{figure}

The cohesive zone, CZ, in a material is the area between two separating but still sufficiently close surfaces ahead of the crack tip, see the shaded region in Fig. \ref{cohesivepic}.

The cohesive forces at the cohesive zones pull the CZ faces together, while
the external load applied to the body, on the contrary, causes  the crack faces and CZ faces to move further apart and the crack to propagate.
When the crack propagates, the cohesive forces vanish at the points where the opening reaches a critical value and these points become the crack surface points, while
the new material points, where the history-dependent normalised equivalent stress reaches a critical value, join the CZ. So, the CZ is practically attached to the crack tip ahead of the crack and moves with the crack, keeping the normalised equivalent stress finite in the body.

One of the most popular CZ models for elasto-perfectly plastic materials is the Leonov-Panasyuk-Dugdale (LPD) model, see \cite{dugdale}, \cite{leo}. In the LPD model, the maximal normal stress in the cohesive zones is constant and equals to the material yield stress, $\sigma=\sigma_y$.
Some generalisations of this model to visco-elastic materials have been developed. In \cite[Section 6.2.3]{Gross-Seelig2011} an approximate solution for the LPD model for a linear visco-elastic bulk material with constant $\sigma=\sigma_y$ on the cohesive zone is presented. Paper \cite{Wnuk-Knauss1970} dealt with penny-shaped crack in a linear visco-elastic bulk material with the ideally rigid-plasticity on the cohesive zone where the yield modulus depends on the CZ tip stain history.   In \cite{Schapery1984}, an energy approach was implemented to the CZ model with constant stress on the cohesive (failure) zone in a non-linear visco-elastic bulk material, while in  \cite{McCartney1988} with an exponential dependence of the CZ stress on the displacement jump, in a linear visco-elastic bulk material. The CZ models with non-constant CZ stress are sometimes associated with the  Barenblatt CZ model, \cite{Barenblatt1959a, Barenblatt1959b, Barenblatt1962}, also described e.g. in \cite[Section 5.3]{Gross-Seelig2011}.

The three main components of the the LPD-type CZ models are: (i) the constitutive equations in the bulk of the material; (ii) the constitutive equations in the CZ; (iii) the criterion for the CZ to break, i.e. for the crack to propagate.

The model presented in this paper is an extension of the LPD model to linear visco-elastic materials with non-linear history-dependent constitutive equations in the CZ.
The history-dependent CZ models can be traced back to  \cite{Wnuk-Knauss1970} but, unlike that paper, the history dependence in our paper is given by a recently developed normalised equivalent stress \cite{MikNam2011AAM} based on the  durability diagram of the material, while the yield condition relates locally with the stress-history of considered CZ points instead of being approximated by the  history at the CZ tip only. The crack starts propagating, when the crack tip opening reaches a critical value (as in the classical LPD model).
Our aim is to find the time evolution of the CZ before the crack starts propagating, the delay time (otherwise known as initiation time), after which the crack will start to propagate, and to study the further time evolution of the crack and the CZ. Some preliminary results on this model were published in \cite{imperial}.

To obtain the solutions, some numerical algorithms were implemented to solve the obtained nonlinear integro-differential problem  and the order of the solution convergence is analysed.

\section{Problem Formulation}

Following \cite{MikNam2011AAM}, let us introduce at any material point the normalised history-dependent equivalent stress
\be
\label{czc1}
\underline{\Lambda}(\boldsymbol{\hat{\sigma}};\hat{t})
=\left(\frac{\beta}{b\sigma_0^{\beta}}\int_0^{\hat{t}} |\boldsymbol{\hat{\sigma}}(\hat{\tau})|^\beta(\hat{t}-\hat{\tau})^{\frac{\beta}{b} -1}d\hat{\tau}\right)^{\frac{1}{\bbb}}
=\left(\frac{\gamma}{\sigma_0^{\beta}}\int_0^{\hat{t}} |\boldsymbol{\hat{\sigma}}(\hat{\tau})|^\beta(\hat{t}-\hat{\tau})^{\gamma-1}d\hat{\tau}\right)^{\frac{1}{\bbb}}
\ee
where
$|\boldsymbol{\hat{\sigma}}|$ is the maximum of the principal stresses, and $\hat{t}$ denotes time.  The parameters $\sigma_0>0$ and $b>0$ are material constants in the assumed power-type relation
\be\nonumber
\hat t_\infty(\hat{\sigma})=\left(\frac{\hat{\sigma}}{\sigma_0}\right)^{-b}
\ee
between the rupture time $\hat t_\infty$ and the constant uniaxial tensile stress $\hat\sigma$ applied to a sample without cracks. The parameter $\beta>0$ is a material constant in the nonlinear accumulation rule for durability under variable load, see \cite{MikNam2011AAM}.
To simplify further formulas, we also introduce the notation $\gamma=\beta/b $.
As shown in \cite{MikNam2011AAM}, the accumulation rule based on \eqref{czc1} becomes equivalent to the Robinson rule of linear summation of the partial life times, when $\beta=b$, i.e., $\gamma=1$.

We will replace the classical LPD CZ (yield) stress condition, $\hat\sigma=\sigma_{\textrm{yield}}$, with the history-dependent condition
\be
\label{czc}
\underline{\Lambda}(\boldsymbol{\hat{\sigma}};\hat{t})=1,
\ee
while in the rest (the bulk) of the material the strength (non-yield) condition should be satisfied
\be
\nonumber
\underline{\Lambda}(\boldsymbol{\hat{\sigma}};\hat{t})<1.
\ee

Note that relations \eqref{czc1},\eqref{czc} were implemented in \cite{MNGeneva03Creep} and \cite{HakMikIMSE10} to solve a similar crack propagation problem without CZ; i.e. it was assumed that when condition (\ref{czc}) is reached at a point, this point becomes part of the crack. However, such approach appeared to be inapplicable for $b\ge 2$, while for many structural materials this parameter is in the range between 5 and 15. In this paper, a CZ approach is developed instead, in order to cover the larger range of $b$ values relevant to structural materials. In the CZ approach, when condition (\ref{czc}) is reached at a point, this point becomes not yet part of  the crack, as in  \cite{MNGeneva03Creep} and \cite{HakMikIMSE10},  but part of the CZ; to became part of the crack, another condition (on the CZ opening) should be satisfied.

As proved further, in Section~\ref{AP-CZ}, the CZ defined by \eqref{czc1}, \eqref{czc} can exist only if $0<\gamma<1$, i.e., $0<\beta<b$. Note that  in all examples considered in \cite{MikNam2011AAM}, where these parameters were calculated by fitting the creep durability experimental data, it was found that $0<\beta<b$, i.e., $0<\gamma<1$. Thus we will further consider $\gamma$ from this interval only.

Let the problem geometry be as in Fig. \ref{cohesivepic}, i.e, the crack occupies the interval $[-\hat{a}(\hat{t}),\hat{a}(\hat{t})]$ and the CZ occupies the intervals $[-\hat{c}(\hat{t}), -\hat{a}(\hat{t})]$ and $[\hat{a}(\hat{t}), \hat{c}(\hat{t})]$ in an infinite linearly elastic or visco-elastic plane loaded at infinity by a traction $\hat q$ in the direction normal to the crack, which is constant in $\hat{x}$,  applied at the time $\hat{t}=0$ and kept constant in time thereafter. The initial CZ is absent, i.e., the CZ tip coordinates coincide with the crack tip coordinates, which are prescribed, $\hat c(0)=\hat a(0)=\hat a_0$, while the functions $\hat c(\hat t)$ and $\hat a(\hat t)$ for time $\hat t>0$ are to be found.

Let $(\hat x, \hat y)$ be coordinates in the Cartesian frame with the origin in the centre of the crack and the $\hat x-$axis directed along the crack.
Assuming that $\hat{\sigma}(\hat{x},\hat{\tau}):=\hat{\sigma}_{\hat y\hat y}(\hat{x},0,\hat{\tau})$ is non-negative and is the maximal component of the stress tensor ahead of the crack,
the CZ condition (\ref{czc})-(\ref{czc1}) at a point $\hat{x}$ on the CZ at time $\hat t$ can be rewritten as
\be
\int_{\hat{t}_c(\hat{x})}^{\hat{t}} \hat{\sigma}^{\beta}(\hat{x},\hat{\tau})(\hat{t}-\hat{\tau})^{\gamma -1} d\hat{\tau}
=\frac{\sigma_0^{\beta}}{\gamma}-\int_0^{\hat{t}_c(\hat{x})} \hat{\sigma}^{\beta}(\hat{x},\hat{\tau})(\hat{t}-\hat{\tau})^{\gamma -1} d\hat{\tau},\label{vol}
\ee
for $\hat{t}\ge \hat{t}_c(\hat{x})$ and $\hat{a}(\hat{t})\le|\hat{x}|\le \hat{c}(\hat{t})$.
Here, $\hat{t}_c(\hat{x})$ denotes the time when the point $\hat{x}$ joined the CZ.
Equation (\ref{vol}) is an inhomogeneous linear Volterra integral equation of the Abel type with unknown function $\hat{\sigma}^\beta(\hat{x},\hat{t})$ for $\hat{t}\ge \hat{t}_c(\hat{x})$.

Let us first consider the case of linear elastic constitutive equations for the bulk of the material.
Applying the results by Muskhelishvili (see \cite{mush}, Section 120), we have for the stresses ahead of the CZ in the elastic material,
\be
\hat{\sigma}(\hat{x},\hat{t})=
\frac{\hat{x}}{\sqrt{\hat{x}^2-\hat{c}^2(\hat{t})}}\left(\hat q -\frac{2}{\pi}\int_{\hat{a}(\hat{t})}^{\hat{c}(\hat{t})} \frac{\sqrt{\hat{c}^2(\hat{t})-\hat{\xi}^2}}{\hat{x}^2-\hat{\xi}^2} \hat{\sigma}(\hat{\xi},\hat{t})d\hat{\xi}\right),
\label{sif}
\ee
for
$|\hat{x}|>\hat{c}(\hat{t})$. As one can see from \eqref{sif}, $\hat{\sigma}(\hat{x},\hat{t})$ has generally a square root singularity as $\hat{x}$ tends to the CZ tip $\hat c$.
The stress intensity factor, $\hat K$, at this singularity can be obtained by multiplying the stress in equation (\ref{sif}) by $\sqrt{\hat{x}-\hat{c}(\hat{t})}$ and taking the limit as $\hat{x}$ tends to $\hat{c}(\hat{t})$, which yields
\be \hat K(\hat{t}) = \sqrt{\frac{\hat{c}(\hat{t})}{2}}\left(\hat q -\frac{2}{\pi} \int_{\hat{a}(\hat{t})}^{\hat{c}(\hat{t})}\frac{\hat{\sigma}(\hat{\xi},\hat{t})}{\sqrt{\hat{c}^2(\hat{t})
-\hat{\xi}^2}}d\hat{\xi}\right). \nonumber\ee

A sufficient condition for the normalised equivalent stress, $\Lambda$, to have no such singularity at the CZ tip is
that the stress $\hat\sigma$ given by \eqref{sif} does not have it either, while the necessary condition for the latter is
that the stress intensity factor, $\hat K$, is zero there.

To simplify the equations, we will employ the following normalisations:
\begin{align}\label{norm}
 & t=\frac{\hat{t}}{\hat t_{\infty}},\quad
 x=\frac{\hat x}{\hat a_0},\quad
 a(t)=\frac{\hat a(t\,\hat t_\infty)}{\hat a_0},\quad c(t)=\frac{\hat{c}(t\,\hat t_\infty)}{\hat a_0},\nonumber\\
 &\sigma(x,t)=\frac{\hat\sigma(x\, \hat a_0,t\, \hat t_\infty)}{\hat q},\quad
 K(c,t)=\frac{\hat K(c\,\hat a_0,t\,\hat t_\infty)}{\hat q\sqrt{\hat a_0}},
\end{align}
where $\hat t_\infty=\hat t_\infty(\hat q)=\left(\frac{\hat{q}}{\sigma_0}\right)^{-b}$ denotes the fracture time for an infinite plane without a crack under the same load, $\hat q$.

Then we obtain the following normalised principle equations for the considered problem:\\
(a) the CZ condition \eqref{czc} in the form
\be
 \int_{t_c(x)}^{t}\sigma^\beta(x,\tau)
(t-\tau)^{\gamma -1} d\tau=\frac{1}{\gamma}
-\int_{0}^{t_c(x)}\sigma^\beta(x,\tau)
(t-\tau)^{\gamma -1}d\tau \
\mbox{ for}\ a(t)\le|x|\le c(t),
\ t>t_c(x);\label{pe1}
\ee
(b) the expression for the stress ahead of the CZ:
\be
\sigma(x,t)=
\frac{x}{\sqrt{x^2-c^2(t)}}
\left(1 -\frac{2}{\pi}\int_{a(t)}^{c(t)} \frac{\sqrt{c^2(t)-\xi^2}}{x^2-\xi^2}\sigma(\xi,t)d\xi\right)
\mbox{ for }|x|>c(t);\label{pe2}
\ee
(c) the zero stress intensity factor, $K(c,t) = 0$ for $t>0$, where
\be
K(c,t) = \sqrt{\frac{c(t)}{2}} -\frac{\sqrt{2c(t)}}{\pi}\int_{a(t)}^{c(t)}\frac{\sigma(\xi,t)}{\sqrt{c^2(t)-\xi^2}}d\xi.\label{pes}
\ee

\section{Cohesive Zone Growth for the Stationary Crack}\label{stat}
In this section we will consider the stationary stage,  when $a(t)=a(0)=1$, and  only the CZ grows with time.  Our aim here is to find the CZ tip position $c(t)$. This stage is followed later by the crack propagation stage considered in the next sections.

\subsection{Numerical Method on the Stationary Crack Stage}\label{nm}
Let us introduce a time mesh with nodes $t_i=ih$, for $i=0,1,2,3,...,n$, where $h=1/n$ is a time increment and $t_n=1$.
At each time step $t_i$, we use the secant method to find the roots, $c(t_i)=c_i$, of the equation $K(c_i,t_i)=0$,
as follows:
\begin{enumerate}
  \item Take 2 initial approximations, $c_{i1}$ and $c_{i2}$, for $c(t_i)$.
  \item Obtain $K_1=K(c_{i1},t_i)$ and $K_2=K(c_{i2},t_i)$ using equation (\ref{pes}).
In order to evaluate the integral in (\ref{pes}),
   \be\label{SIFI}
 \int_{a(0)}^{c(t_i)}\frac{1}{\sqrt{c^2(t_i)-\xi^2}}\sigma(\xi,t_i)d\xi,
 \ee
   we piece-wise linearly interpolate $\sigma(\xi,t_i)$ in $\xi$ on the CZ, over the points $c(t_k)$, $k=0,1,2,...,i$.

To obtain $\sigma(c(t_i),t_i)$, we use the Abel integral equation (\ref{pe1}) with zero left hand side, which reduces the equation to the following one,
    $$\int_{0}^{t_i}{\sigma^\beta(c(t_i),\tau)}
(t_i-\tau)^{\gamma -1}d\tau=\frac{1}{\gamma}.$$
By explicit integration of the
   piece-wise linear interpolant of the function $\sigma^\beta(c(t_i),\tau)$ in $\tau$ over the time instants $t_j$, $j=0,1,2,...,i$, including the unknown value $\sigma^\beta(c(t_i),t_i)$, we obtain the linear algebraic equation for the latter, which has the solution
\begin{multline}
\sigma^\beta(c(t_i),t_i)=
\sigma^\beta(c(t_i),t_{i-1})+
(t_i-t_{i-1})^{-\gamma}\left\{\vphantom{\sum_{j=0}^{i-2}}
(\gamma+1)[1-\sigma^\beta(c(t_i),0)t_i^\gamma]\ +\right.\\
\sum_{j=0}^{i-2}\sigma^\beta(c(t_i),t_j)\frac{(t_i-t_j)^{\gamma+1}-(t_i-t_{j+1})^{\gamma+1}}{t_{j+1}-t_j}\
+\\
\left.\sum_{j=1}^{i-1}\sigma^\beta(c(t_i),t_j)\frac{(t_i-t_j)^{\gamma+1}-(t_i-t_{j-1})^{\gamma+1}}{t_j-t_{j-1}}
\right\}.
\label{sti}
\end{multline}

To obtain $\sigma(c(t_k),t_i)$, at each $k<i$, we use the Abel integral equation (\ref{pe1}) with $x=c(t_k)$, $t_c(x)=t_k$ and $t=t_i>t_k$.
First, we again evaluate the right hand side integral
    $$\int_{0}^{t_k}{\sigma^\beta(x,\tau)}
(t-\tau)^{\gamma -1}d\tau$$
by analytic integration of the piece-wise linear interpolant of the function $\sigma^\beta(x,\tau)$ in $\tau$ over the time instants $t_j$, $j=0,1,2,...,k$, where $t_0=0$.
Then we use the analytical solution (\ref{sol2}) of the Abel-type integral equation (\ref{pe1}), see details in Subsection \ref{Sab}, and arrive at the following solution at $t>t_k$ for $x=c(t_k)$,
\begin{multline}\sigma^{\beta}(x,t)=-\frac{1}{\pi}\sin(\pi\gamma) \sum_{j=1}^{k}\left\{\sigma^{\beta}(x,t_{j-1})\left[V(t_{j-1},t,t_k)-V(t_j,t,t_k)\right]+
\right. 
\\
\left.\frac{\sigma^{\beta}(x,t_j)-\sigma^{\beta}(x,t_{j-1})}{\gamma(t_j-t_{j-1})} \left[W(t_{j-1},t,t_k)-W(t_j,t,t_k)-\gamma (t_j-t_{j-1})V(t_j,t,t_k)\right] \right\}\\
\shoveleft{=-\frac{1}{\pi}\sin(\pi\gamma) \sum_{j=1}^{k}\left\{\sigma^{\beta}(x,t_{j-1})
\left[\widetilde V(t_{j-1},t,t_k)-\widetilde V(t_j,t,t_k)\right]+
\right. 
}\\
\left.\frac{\sigma^{\beta}(x,t_j)-\sigma^{\beta}(x,t_{j-1})}{\gamma(t_j-t_{j-1})}
\left[\widetilde W(t_{j-1},t,t_k)-\widetilde W(t_j,t,t_k)-\gamma (t_j-t_{j-1})\widetilde V(t_j,t,t_k)\right] \right\}\\
\shoveleft{=-\text{sinc}(\pi\gamma)\left\{
\sigma^{\beta}(x,0)\left[\widetilde{\widetilde V}_0(t,t_k)
+\frac{\widetilde W(t_{1},t,t_k)-\widetilde W(0,t,t_k)}{t_1}\right]\right.
}\\
+\sum_{j=1}^{k-1}\sigma^{\beta}(x,t_{j})
\left[\frac{\widetilde W(t_{j+1},t,t_k)-\widetilde W(t_j,t,t_k)}{t_{j+1}-t_j}
-\frac{\widetilde W(t_{j},t,t_k)-\widetilde W(t_{j-1},t,t_k)}{t_{j}-t_{j-1}}\right]\\
\left. + \sigma^{\beta}(x,t_{k})\frac{\widetilde W(t_{k-1},t,t_k)}{t_{k}-t_{k-1}}
\right\}\\
\shoveleft{=-\text{sinc}(\pi\gamma)\left\{
\sigma^{\beta}(x,0)\widetilde{\widetilde V}_0(t,t_k)
+\frac{\sigma^{\beta}(x,t_1)-\sigma^{\beta}(x,0)}{t_1}\widetilde W(0,t,t_k) \vphantom{\sum_{j=1}^{k-1}}\right.
}\\
\left.+\sum_{j=1}^{k-1}
\left[\frac{\sigma^{\beta}(x,t_{j+1})-\sigma^{\beta}(x,t_{j})}{t_{j+1}-t_{j}}
-\frac{\sigma^{\beta}(x,t_{j})-\sigma^{\beta}(x,t_{j-1})}{t_{j}-t_{j-1}}\right]\widetilde W(t_{j},t,t_k)
\right\},
\label{ab}
\end{multline}
where
\begin{align}\label{Ve}
&V(y,t,t_c)= \int_{t_c}^{t}\frac{(\tau-y)^{\gamma -1}}{(t-\tau)^{\gamma }}d\tau
=\pi\csc{\left(\pi\gamma \right)} +\widetilde V(y,t,t_c), \\
\label{Vet}
&\widetilde V(y,t,t_c)=
-\frac{1}{\gamma}\left(\frac{t_c-y}{t-y}\right)^{\gamma }\,\, _2F_1 \left[\gamma ,\gamma ;1+\gamma ;\frac{t_c-y}{t-y}\right],\\
&\widetilde{\widetilde V}_0(t,t_k)=\gamma\widetilde V(0,t,t_c)=
-\left(\frac{t_c}{t}\right)^{\gamma }\,\, _2F_1 \left[\gamma ,\gamma ;1+\gamma ;\frac{t_c}{t}\right],\\
\label{We}
& W(y,t,t_c)=\int_{t_c}^{t}\frac{(\tau-y)^{\gamma}}{(t-\tau)^{\gamma }}d\tau
=\gamma \pi\csc{\left(\pi\gamma \right)}(t-y)+\widetilde W(y,t,t_c),\\
&\widetilde W(y,t,t_c)=- \frac{1}{1+\gamma}(t_c-y)^{1+\gamma }(t-y)^{-\gamma } \,\,_2F_1 \left[1+\gamma ,\gamma ;2+\gamma ;\frac{t_c-y}{t-y}\right],
\end{align}
$_2F_1$ is the Gauss hypergeometric function, and $\text{sinc}(\pi\gamma)=\dfrac{\sin(\pi\gamma)}{\pi\gamma}$.

To implement \eqref{sti} and \eqref{ab}, we need, in turn, to find $\sigma(c(t_k),t_j)$ for $0\le t_j<t_k\le t_i$ from equation (\ref{pe2}) (since $c(t_k)>c(t_j)$), $j=0,1,\dots,k-1$. For $j=0$, $t_0=0$ and the integral in (\ref{pe2}) vanishes giving $\sigma(c(t_k),0)=c(t_k)/\sqrt{c^2(t_k)-1}$.
For $j>0$, taking into account that $K(c(t_j),t_j)=0$, equation (\ref{pe2}) reduces to
\be
\sigma(x,t)=
\frac{2x}{\pi}\sqrt{x^2-c^2(t_j)}\int_{a(t)}^{c(t_j)} \frac{\sigma(\xi,t_j)}{(x^2-\xi^2)\sqrt{c^2(t_j)-\xi^2}}d\xi\
\mbox{ for }|x|>c(t_j),\label{pe2r}
\ee
where the integral is calculated,  similarly to  integral \eqref{SIFI},  linearly interpolating $\sigma(\xi,t_j)$ between $\xi=c(t_m)$ and $\xi=c(t_{m+1})$ for $m=0,1,...,j-1$.


  \item Find the next approximation for $c_i$ using
  $$(c_i)_3=\frac{K_2 \cdot c_{i1}-K_1\cdot c_{i2}}{K_2-K_1}$$
  \item If $|(c_i)_3-c_{i1}|<\epsilon$ or $|(c_i)_3-c_{i2}|<\epsilon$ then convergence is reached and we allocate $c(t_i)=c_3$ and go to the step $t=t_{i+1}$; otherwise, take the new $c_{i2}$ as $(c_i)_3$ and return to item 2. Here $\varepsilon$ is some tolerance.
\end{enumerate}

\subsection{Numerical Results for the Stationary Crack}
Programming of the described algorithm was implemented in MATLAB with $\varepsilon=10^{-8}$ as the tolerance value.

\begin{figure}[H]
\centering
\mbox{\subfigure{\includegraphics[scale=0.55]{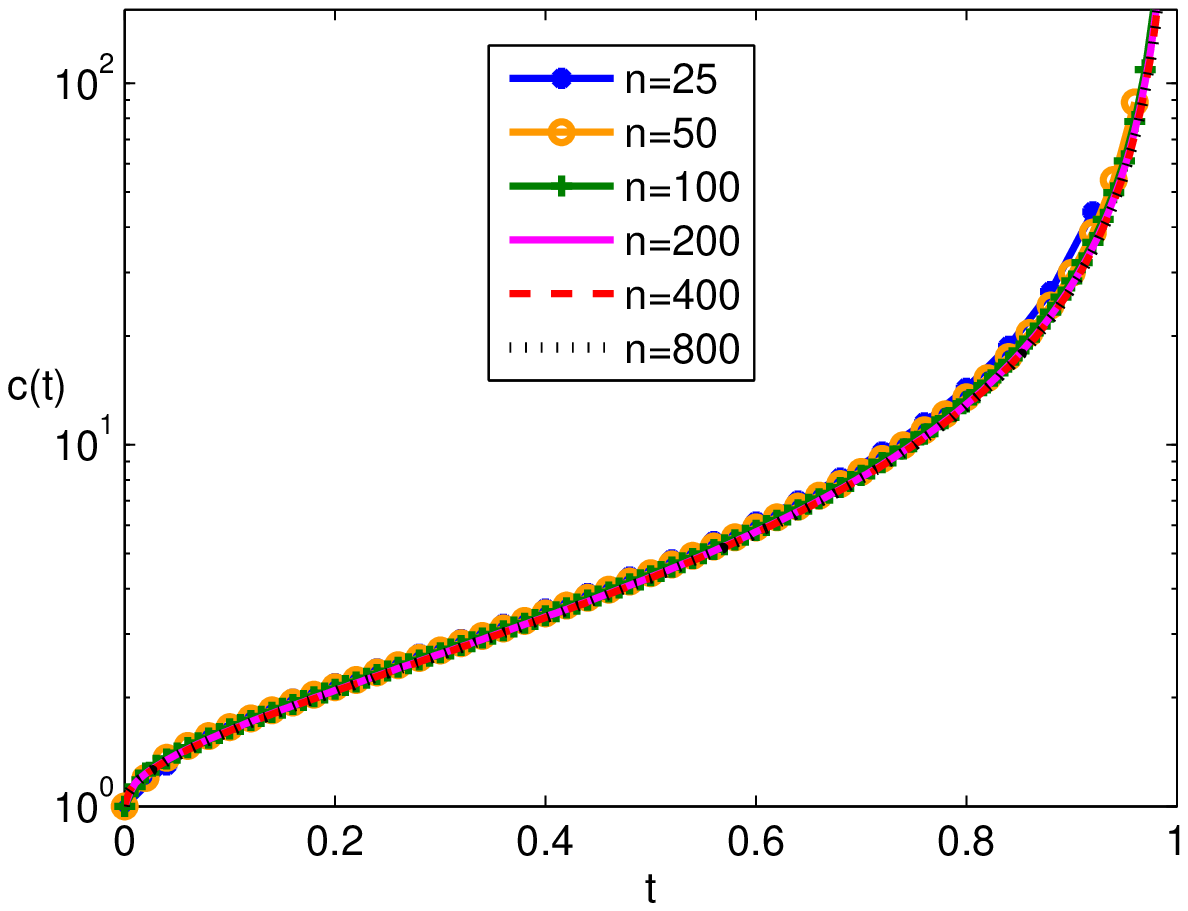}}\,\,\,
\subfigure{\includegraphics[scale=0.55]{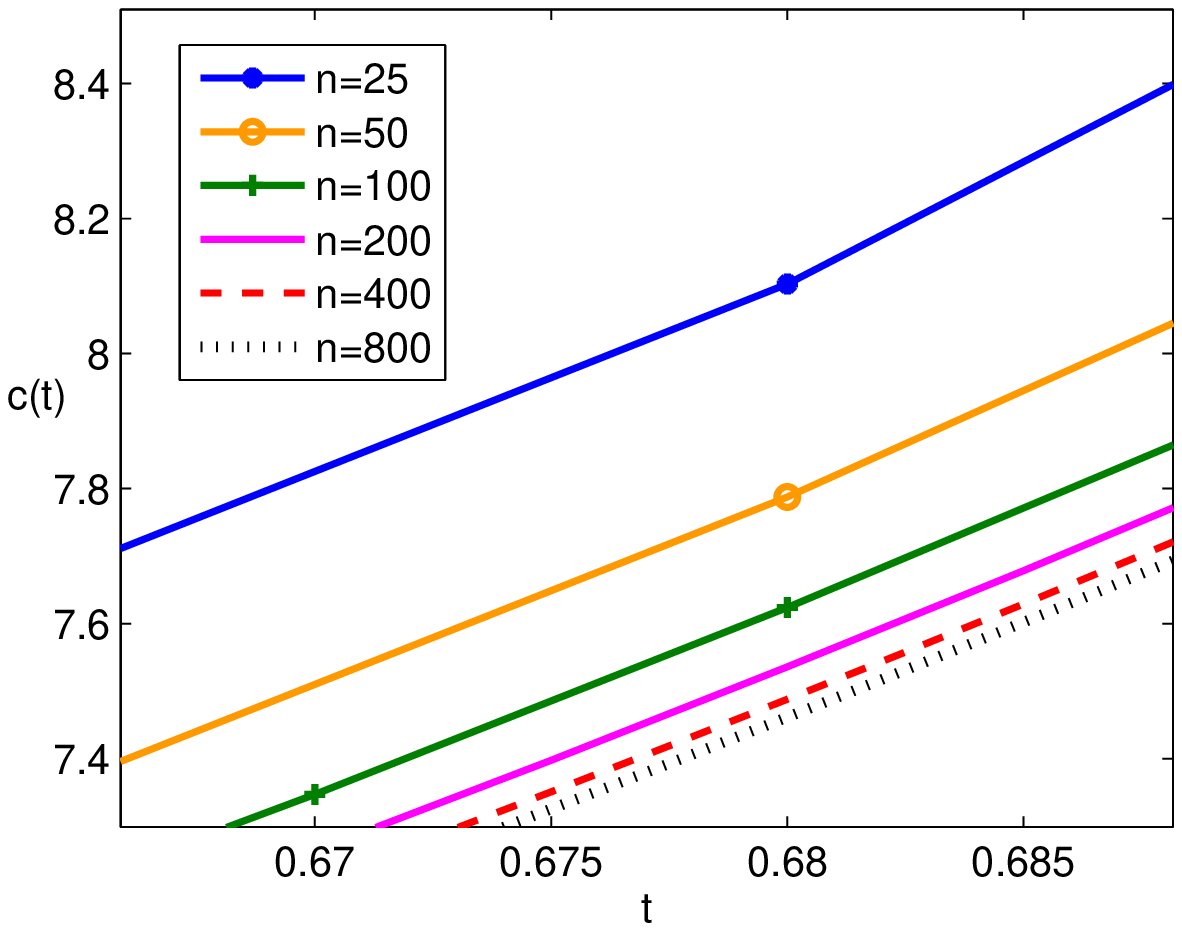}}}
\caption{CZ tip position vs time for $b=4$, $\beta=2$ and different meshes.
}\label{cvt1}
\end{figure}

\begin{figure}[H]
\centering
\begin{minipage}{.45\textwidth}
  \centering
  \includegraphics[width=1\linewidth]{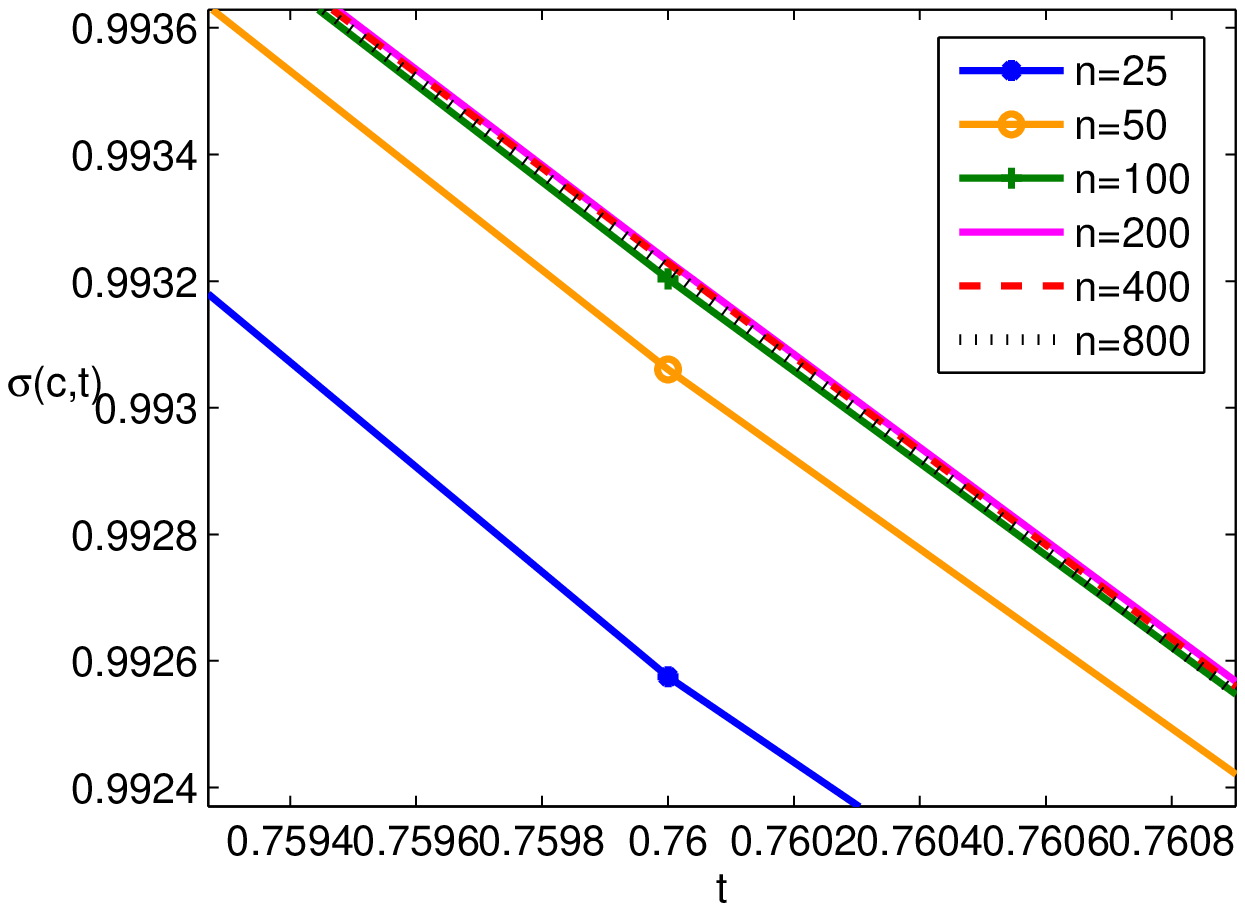}
  \caption{$\sigma(c(t^*),t)$  for $b=4$, $\beta=2$, $t^*=0.6$: closer look in the CZ ($t>t^*$).}
  \label{fig:test1}
\end{minipage}
\quad
\begin{minipage}{.45\textwidth}
  \centering
  \includegraphics[width=1\linewidth]{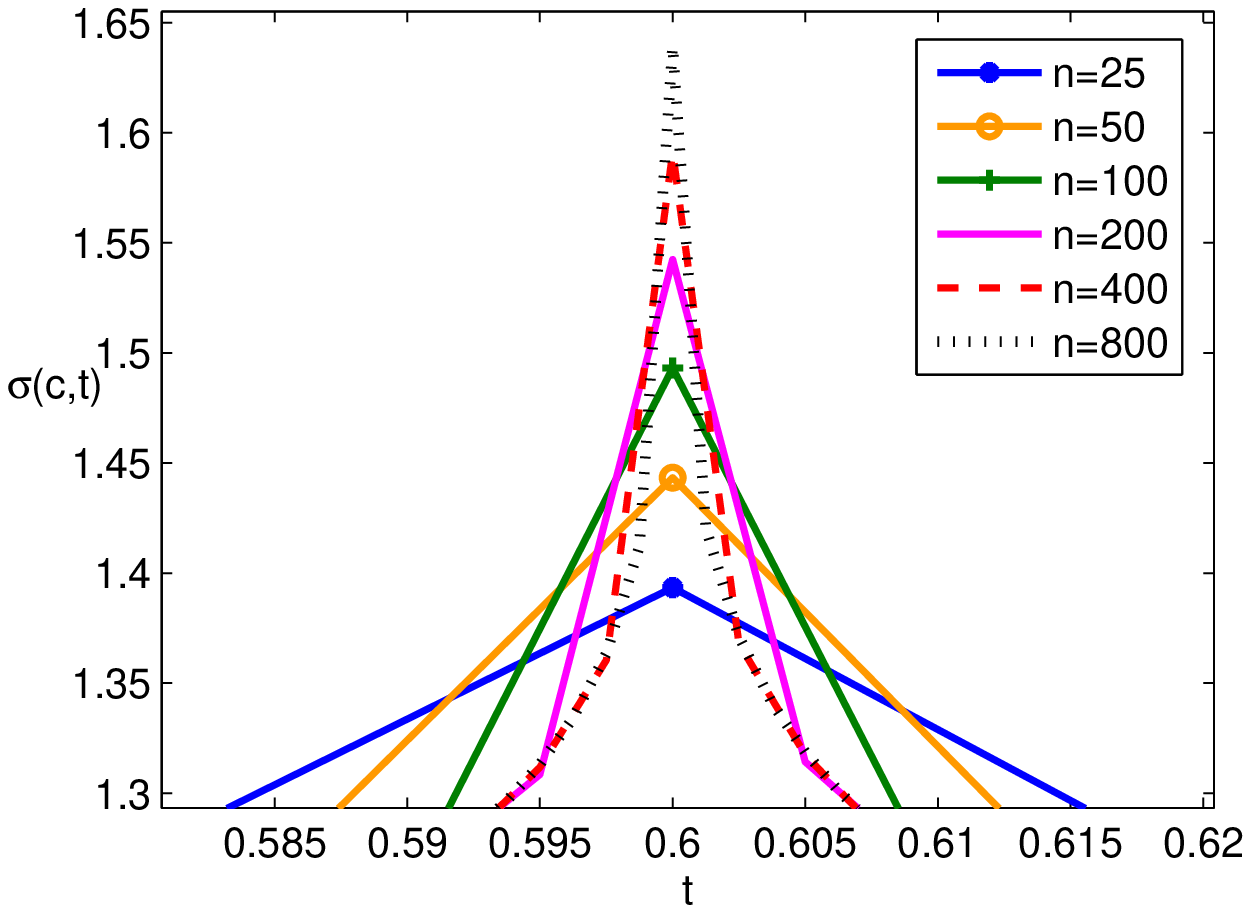}
  \caption{$\sigma(c(t^*),t)$  for $b=4$, $\beta=2$, $t^*=0.6$: closer look at the CZ tip (near $t^*$).}
  \label{fig:test2}
\end{minipage}
\end{figure}

The graphs on Figs.~\ref{cvt1}-\ref{fig:test2} show the obtained numerical results on the evolution of the CZ tip position as well as the stress distribution on the CZ  for various mesh sizes.

\begin{figure}[H]
\centering
\begin{minipage}{.45\textwidth}
  \centering
  \includegraphics[width=1\linewidth]{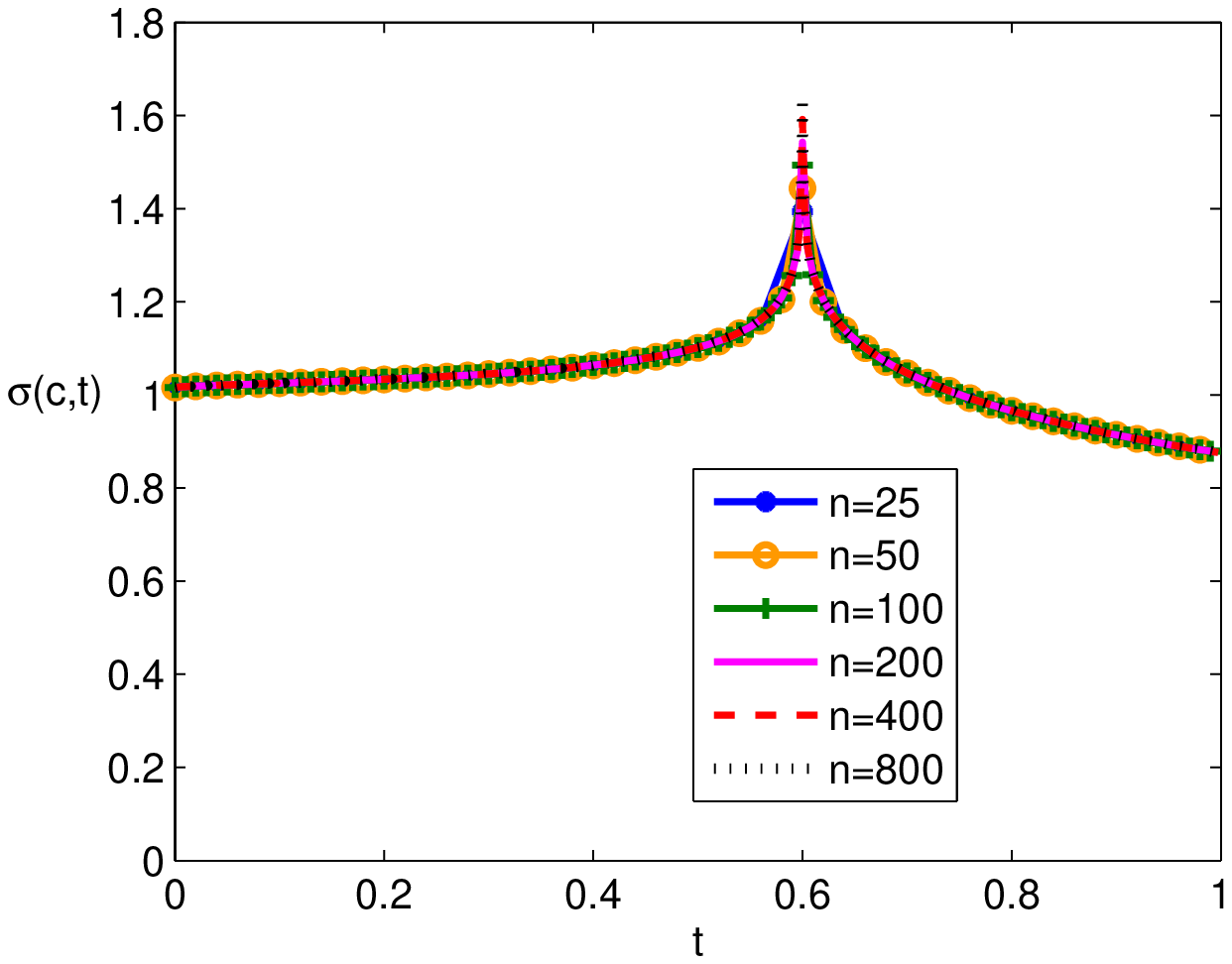}
  \caption{$\sigma(c(t^*),t)$ vs time for $b=4$, $\beta=2$, $t^*=0.6$: global picture.
    }
\end{minipage}
\quad
\begin{minipage}{.45\textwidth}
  \centering
  \includegraphics[width=1\linewidth]{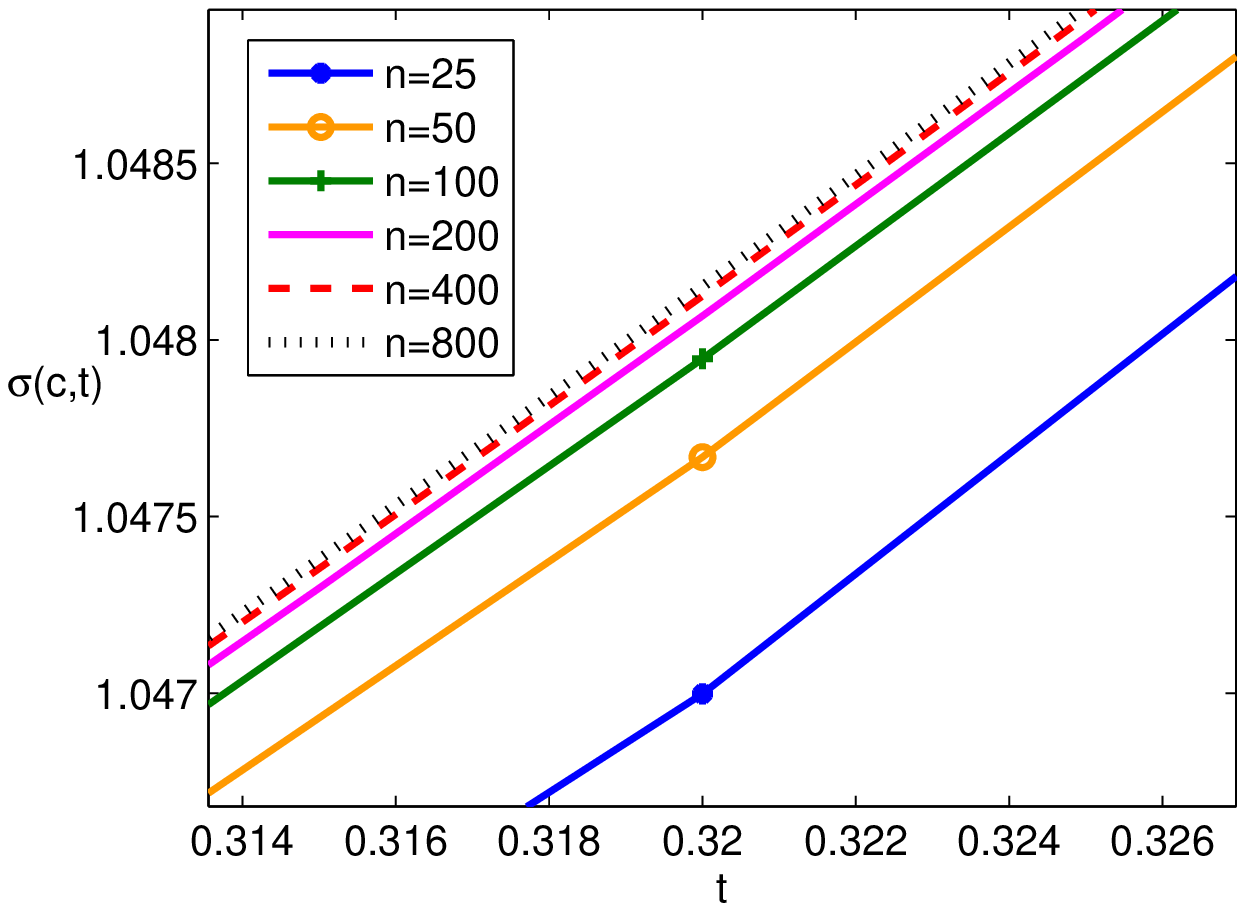}
  \caption{$\sigma(c(t^*),t)$  for $b=4$, $\beta=2$, $t^*=0.6$: closer look ahead of the CZ ($t<t^*$).}
\end{minipage}
\end{figure}

All the graphs illustrate fast numerical convergence of the obtained results, except the graphs in Fig.~\ref{fig:test2} for the stress $\sigma$ at the cohesive zone tip. In more details the convergence in these and the following graphs is analysed in Section~\ref{Conv}.

\section{Crack Tip Opening}
We will first consider the case when the  bulk of the material is linearly elastic and then convert the obtained solution to the case of linear visco-elastic materials using the so-called Volterra principle.
Using the representations by Muskhelishvili (see \cite{mush}, Section 120), it can be deduced that in the linearly elastic isotropic homogeneous plane with a crack,  the normal displacement jump at a point $\hat x$ on the crack or CZ shore  is
\be
[\hat u_e](\hat{x},\hat{t}) =[\hat u_e^{(q)}](\hat{x},\hat{t})
+[\hat u_e^{({\sigma})}] (\hat{x},\hat{t}), \quad |\hat x|<\hat c(t), \nonumber
\ee
where
\be
[\hat u_e^{(q)}](\hat{x},\hat t)= \frac{\hat q(1+\varkappa)}{2\mu_0}\sqrt{\hat{c}(\hat t)^2-\hat{x}^2},\quad
[\hat u_e^{({\sigma})}](\hat{x},\hat{t})= \frac{1+\varkappa}{2\pi \mu_0}\left(\int_{\hat{a}(\hat t)}^{\hat{c}(\hat t)} \hat\sigma(\hat{\xi},\hat{t})\Gamma(\hat{x},\hat{\xi};\hat{c}(\hat t))d\hat{\xi} \right), \nonumber
\ee
and
\be
\Gamma(\hat{x},\hat{\xi};\hat{c})= \ln{\left[\frac{2\hat{c}^2-\hat{\xi}^2-\hat{x}^2 -2\sqrt{(\hat{c}^2-\hat{x}^2) (\hat{c}^2-\hat{\xi}^2)}}{2\hat{c}^2-\hat{\xi}^2-\hat{x}^2+ 2\sqrt{(\hat{c}^2-\hat{x}^2)(\hat{c}^2-\hat{\xi}^2)}}\right]}. \nonumber
\ee
In the above expressions, $\varkappa=3-4\nu$ under the plain strain condition, while $\varkappa=(3-\nu)/(1+\nu)$ under the plain stress condition, $\mu_0=E_0/[2(1+\nu)]$ is the shear modulus, where $E_0$ and $\nu$ denote Young's modulus of elasticity and Poisson's ratio, respectively.


Then the displacement jump  at the crack tip,
that we call the crack tip opening, for elastic material is given by the formula
\be\label{de}
\hat{\delta}_e(\hat t):=[\hat u_e](\hat{a}(\hat t),\hat{t})
=\frac{1+\varkappa}{2 \mu_0}\left(\hat q\sqrt{\hat{c}^2(\hat{t})-\hat{a}^2(\hat{t})}
+\frac{1}{\pi}\int_{\hat{a}(\hat{t})}^{\hat{c}(\hat{t})}\hat{\sigma}(\hat{\xi},\hat{t}) \Gamma(\hat{a}(\hat{t}),\hat{\xi},\hat{c}(\hat{t}))d\hat{\xi}\right).
\ee
Using the space and time normalisations given in equation (\ref{norm}) as well as the following normalisation

\be\label{den}
[u_e(x,t)]=\frac{2 \mu_0}{(1+\varkappa)\hat q}\ \frac{[\hat u_e]\left(\hat x\,\hat a_0,t\hat t_\infty\right)}{\hat a_0},\quad
\delta_e(t)=\frac{2 \mu_0}{(1+\varkappa)\hat q}\ \frac{\hat{\delta}_e(t\, \hat t_\infty)}{\hat a_0},
\ee
we obtain
\begin{align}
\nonumber
&[u_e](x,t)=  \sqrt{c^2(t)-x^2} +
\frac{1}{\pi}\int_{a(t)}^{c(t)}\sigma(\xi,t) \Gamma(x,\xi;c(t))d\xi,\\
\delta_e(t)= &[u_e](a(t),t)= \sqrt{c^2(t)-a^2(t)} +
\frac{1}{\pi}\int_{a(t)}^{c(t)}\sigma(\xi,t) \Gamma(a(t),\xi;c(t))d\xi.\label{deltaelas}
\end{align}

To obtain the crack tip opening in the visco-elastic case, we will implement the so-called Volterra principle, according to which we have to replace  the elastic constants $\mu_0$ and $\nu$ in the elastic solution by the corresponding visco-elastic operators, to arrive at the visco-elastic solution. Although this approach does not always bring a visco-elastic solution for the problems with moving boundaries, it is possible to show, cf. \cite{robot}, that this approach leads to a visco-elastic solution for the plane symmetric problem with a straight propagating crack. This means that for the visco-elastic problem we can directly use the results by Muskhelishvili for the stress representation given in equation (\ref{sif}) since they do not include the elastic constants at all.

For simplicity, we will consider the visco-elastic material with constant (purely elastic) Poisson's ratio $\nu$ (and thus the parameter $\varkappa$). Then, to obtain the crack opening in the visco-elastic case, we have to replace ${1}/{\mu_0}$ in \eqref{de} by the second kind Volterra integral operator $\boldsymbol{\mu}^{-1}$ defined as
\be
\left(\boldsymbol{\mu}^{-1}\hat{\sigma}\right)\left(\hat{t}\right)= \frac{1}{\mu_0}\left\{\hat{\sigma}\left(\hat{t}\right)+ \int_0^{\hat{t}}\dot{\mathcal J}\left(\hat{t}-\hat{\tau}\right) \sigma\left(\hat{\tau}\right)d\hat{\tau}\right\},\nonumber
\ee
where the creep function $\mathcal J$ is known and $\dot{\mathcal J}$ is its derivative, while $\mu_0$ is the instant shear modulus. Hence the visco-elastic crack tip opening becomes
\be
\hat{\delta}_v\left(\hat{t}\right)=[\hat u_v](\hat{a}(\hat t),\hat t)
= \left(\boldsymbol{\mu}^{-1}\mu_0[\hat u_e](\hat{a}(\hat t),\cdot)\right) \left(\hat{t}\right)
= \hat{\delta}_e\left(\hat{t}\right)+ \int_0^{\hat{t}} \dot{\mathcal J}\left(\hat{t}-\hat{\tau}\right)
[\hat u_e](\hat{a}(\hat t),\hat\tau)d\hat{\tau}. \label{3}
\ee

In our numerical examples we use the creep function of a standard linear solid,
\mycomment{
\be \label{JKV}
\mathcal J\left(\hat{t}-\hat{\tau}\right)
= 1+\frac{\mu_0}{\mu_1} \left(1-e^{-\frac{\mu_1}{\eta}\left(\hat{t}-\hat{\tau}\right)}\right).
 \ee
 or
 } 
 \be \label{JKV}
  \dot{\mathcal J}\left(\hat{t}-\hat{\tau}\right)
 = \frac{\mu_0}{\eta} e^{-\frac{\hat t-\hat\tau}{\hat\theta}}.
  \ee
Here,
the material parameters $\hat\theta$ and $\eta$
are, respectively, the relaxation time and the viscosity of the visco-elastic material.
Such visco-elastic models satisfactorily describe some polymers, e.g. PMMA (also known as plexiglas).
For $\dot{\mathcal J}$ in the form \eqref{JKV}, after employing the normalised parameters
\be
\delta_v(t)= \frac{2 \mu_0\hat{\delta}_v(t\, \hat t_\infty)}{(1+\varkappa)\hat a_0 \hat q},\quad  \theta=\frac{\hat\theta}{\hat t_\infty},
\quad m=\frac{\mu_0\hat t_\infty}{\eta},
\label{A0A1}
\ee
equation \eqref{3} reduces to the following expression for the normalised crack tip opening in the visco-elastic case,
\be\label{d=dcvisc2}
\delta_v(t)=[u_v({a}(t),t)]=
\left(\delta_e(t)+m\int_{t_c(a(t))}^{t} e^{-\frac{t-\tau}{\theta}}[u_e]({a}(t),{\tau})d\tau\right),
\ee
where the lower limit of the integral is replaced with $t_c(a(t))$ since $[u_e](x,{\tau})=0$ when $\tau\le t_c(x)$.

In the numerical examples for the visco-elastic case we used values  $m=5$ and $\theta=1$, which are of the order of the ones for PMMA, see Section~\ref{RdPMMA}.

\begin{figure}[H]
\centering
\mbox{\subfigure{\includegraphics[scale=0.5]{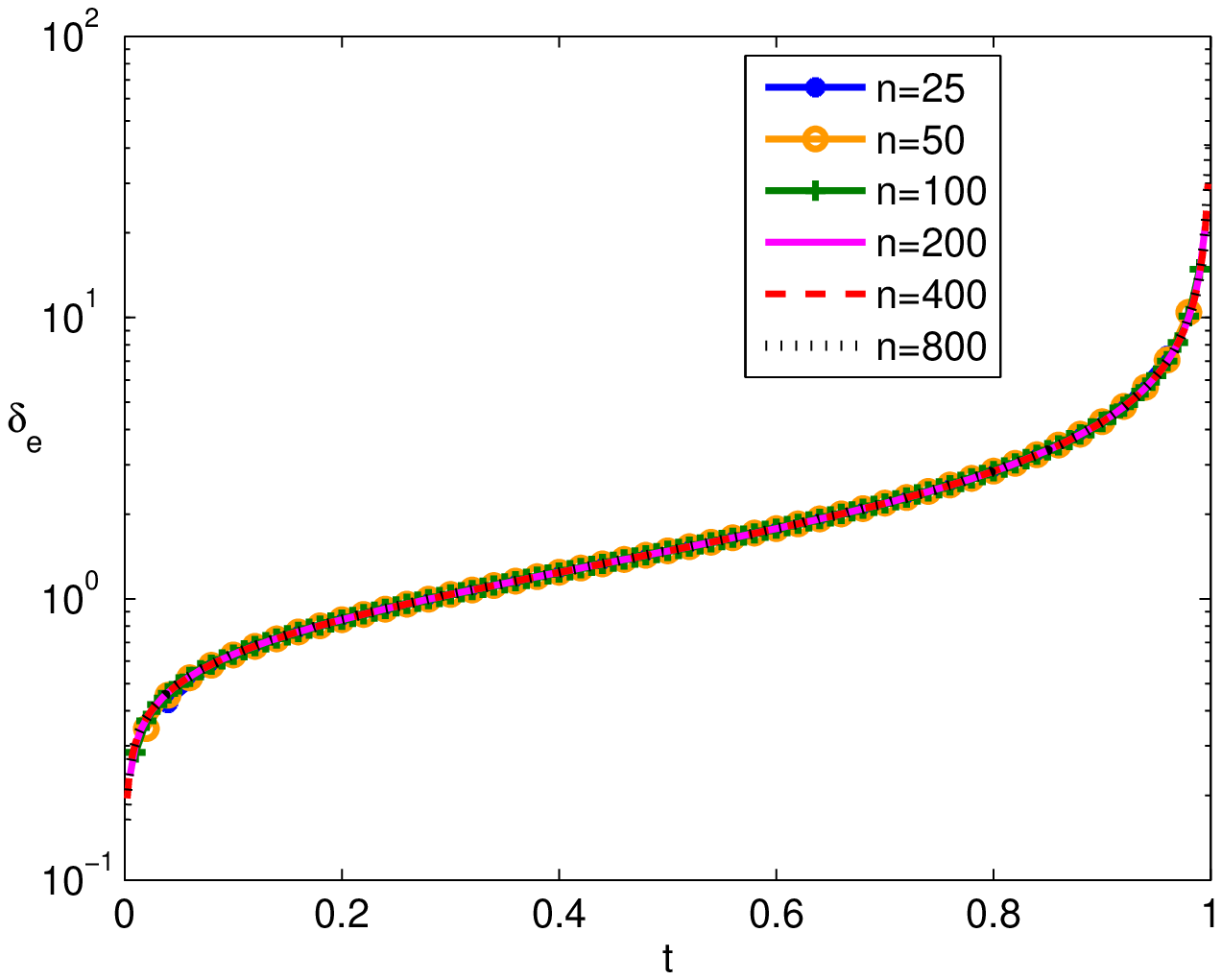}}\quad
\subfigure{\includegraphics[scale=0.5]{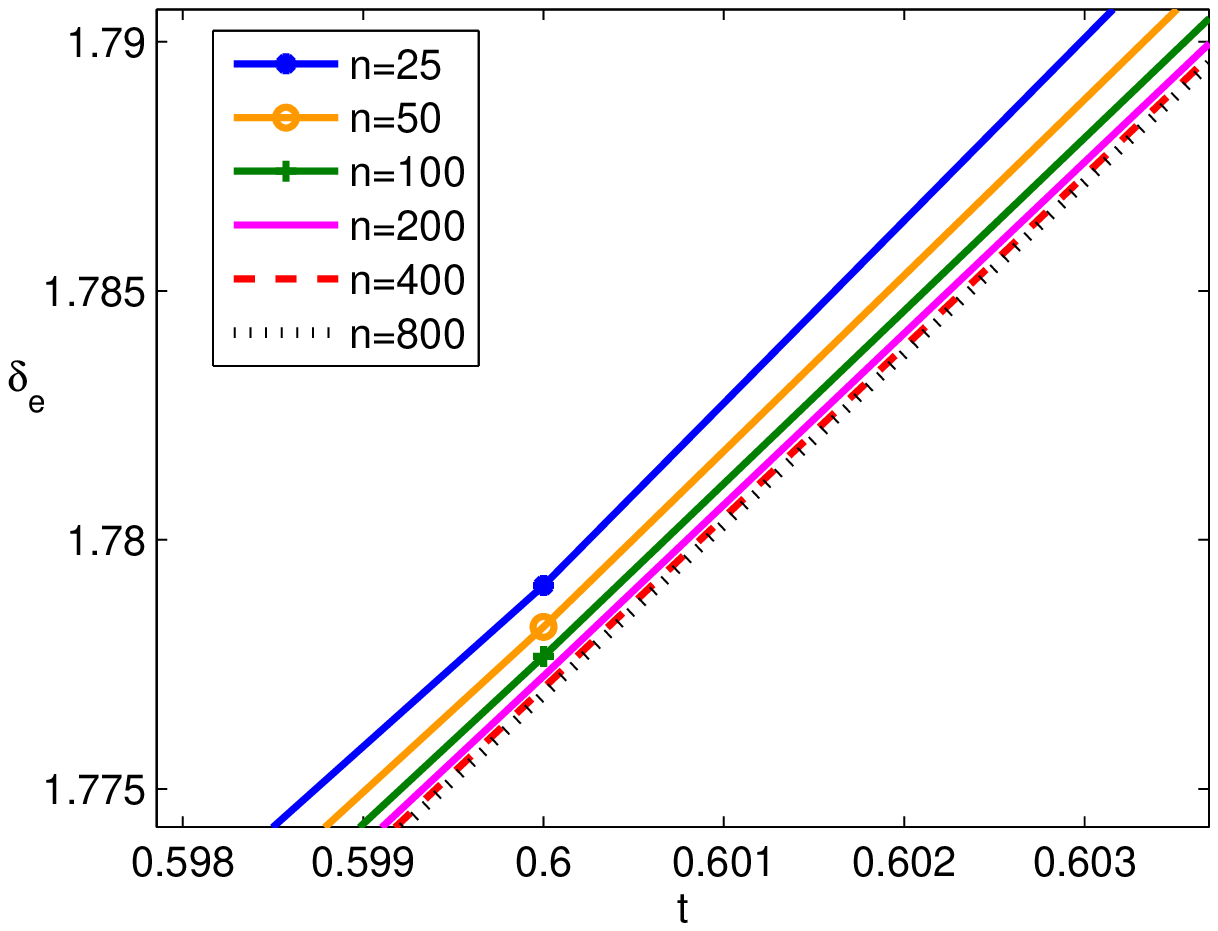}}}
\caption{Crack tip opening $\delta_e$ vs. time $t$ in the elastic case for different time meshes
}\label{co}
\end{figure}

\begin{figure}[H]
\centering
\mbox{\subfigure{\includegraphics[scale=0.5]{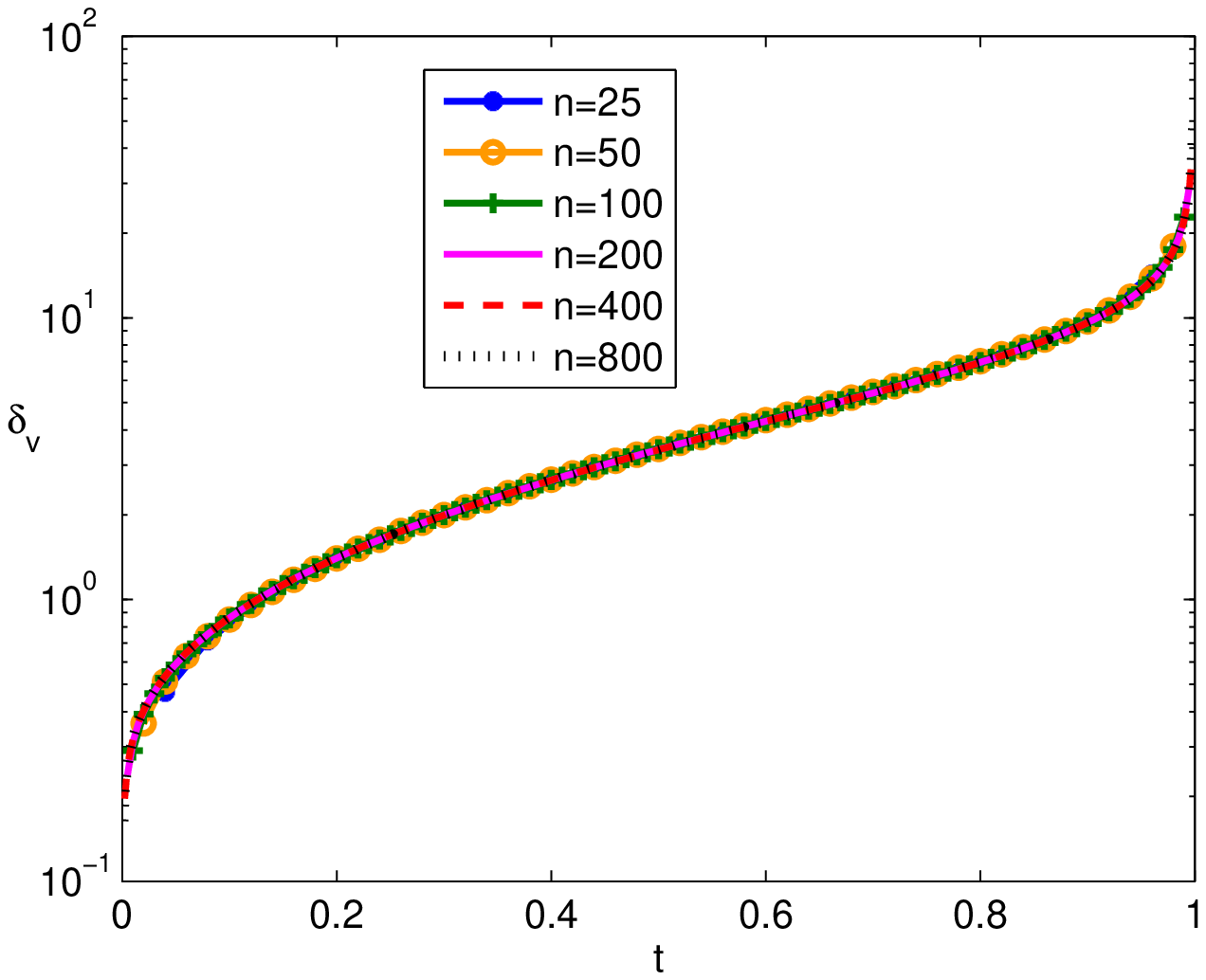}}\quad
\subfigure{\includegraphics[scale=0.5]{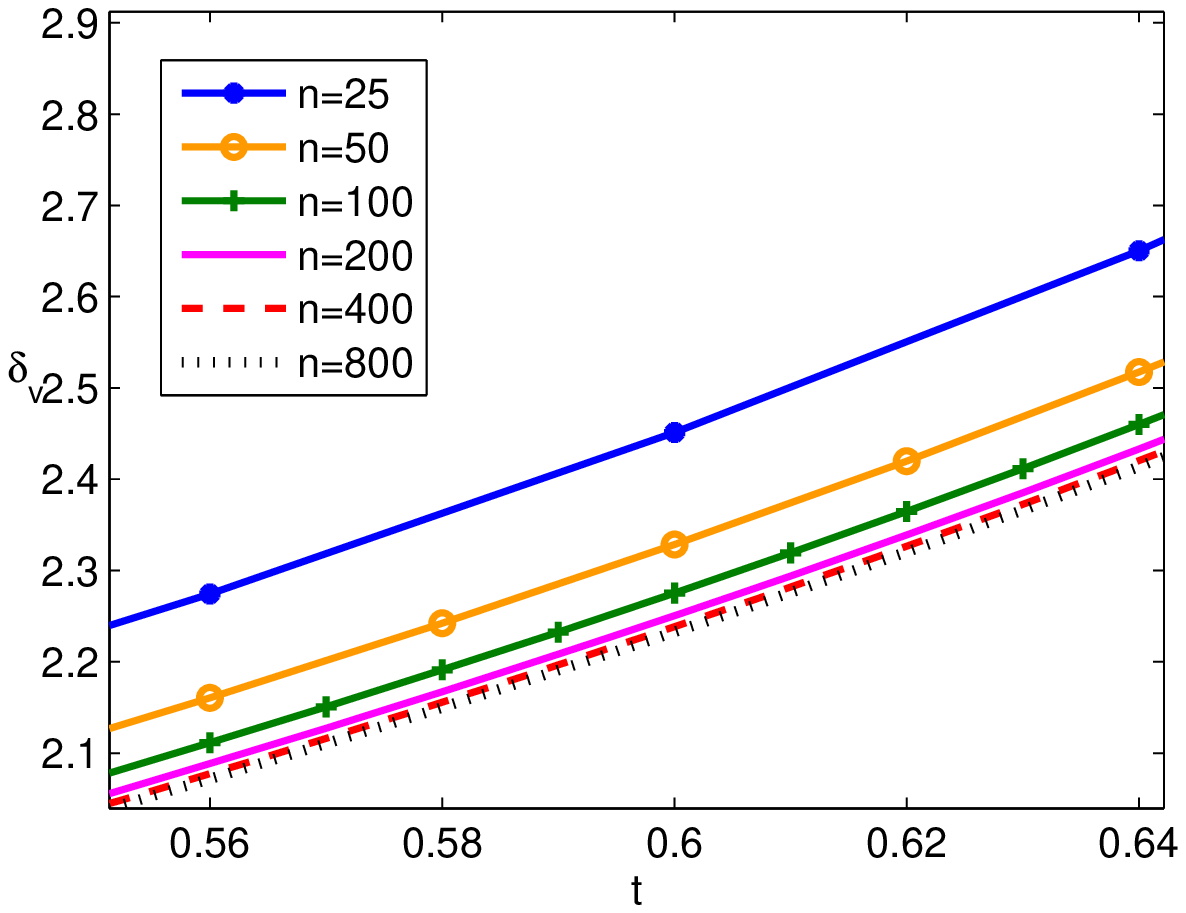}}}
\caption{Crack tip opening $\delta_v$ vs. time $t$ in the visco-elastic case for different time meshes
}\label{co2}
\end{figure}

\begin{figure}[H]
\begin{center}
\includegraphics[scale=0.55]{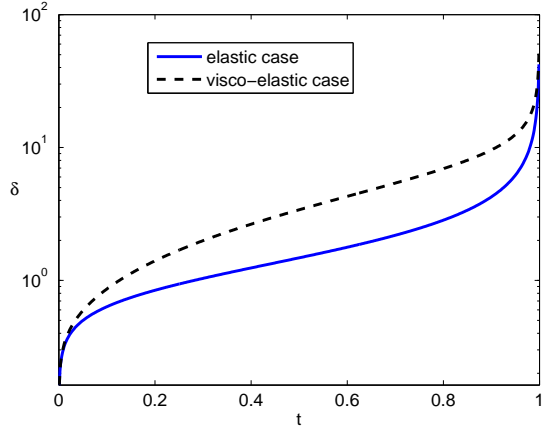}
\caption{Comparison of the crack tip opening $\delta$ vs. time $t$ for elastic and visco-elastic cases}
\label{con}
\end{center}
\end{figure}

The graphs in Figs. \ref{co}-\ref{co2} show the stationary crack tip opening evolution for $b=4$ and $\beta=2$, in the elastic and visco-elastic cases for different time meshes, while Fig.~\ref{con} gives their comparison for the finest mesh, $n=800$.
\section{Crack Propagation Stage}
We have, so far, assumed that the crack is stationary and only the CZ is growing ahead of the crack.
However, the crack will start to propagate when the crack tip opening $\hat{\delta}$ reaches a critical value $\hat{\delta}_c$, which is considered as a material constant.
The crack and CZ tip will not necessarily propagate at the same rate, i.e., the crack length can vary.

Assuming that the external load $\hat q$ is applied at time $\hat t=0$, the time instant, when the crack tip opening reaches a critical value and the crack starts propagating, will be referred to as the fracture delay time (sometimes also named as the fracture initiation time), $\hat t_d$.

Similar to \eqref{norm}, \eqref{den} and \eqref{A0A1}, we employ the following normalised parameters,
\be\label{dc}
t_d=\frac{\hat t_d}{\hat t_\infty},\quad
\delta_c=\frac{2 \mu_0}{(1+\varkappa)\hat q}\ \frac{\hat{\delta}_c}{\hat a_0 }.
\ee

The crack tip opening $\delta(t_i)$ satisfies equation
\be \label{d=dc}
\delta_e(t)=\delta_c, \quad t\ge t_d.
\ee
for the purely elastic case or equation
\be \label{d=dcvisc}
\delta_v(t)=\delta_c, \quad t\ge t_d.
\ee
for the visco-elastic case, where $\delta_e(t_i)$ and $\delta_v(t_i)$ are given by \eqref{deltaelas} and \eqref{d=dcvisc2}, respectively.

\subsection{Numerical Algorithm on the Propagating Crack Stage}\label{RDNAPCS}
The first aim in the crack propagation stage is to find the delay time, $t_d$, solving equations (\ref{d=dc}) or (\ref{d=dcvisc}) by the secant method, while $a(t_d)=1$ and the corresponding value of $c(t_d)$ is obtained by setting the stress intensity factor to zero and applying secant iterations as explained in Section \ref{nm}.
Then, to calculate the crack length and the CZ length at $t>t_d$, we
use the uniform time mesh with time steps $t_i=t_d+i\cdot h$, where $h$ is the step size,
and implement the secant method to solve equation \eqref{d=dc} (in the elastic case) or \eqref{d=dcvisc} (in the visco-elastic case) for $a(t_i)$. To do this, we need  $c(t_i)$ at each iteration, which is obtained using the
secant method to solve equation $K(c(t_i),t_i)=0$
for $c(t_i)$, where the stress intensity factor $K(c,t)$ is given by (\ref{pes}). Further details on the algorithm are given below.
\subsubsection{Initial Approximations}
Note that we take the following 2 initial approximations for $a_i:=a(t_i)$: $(a_i)_1=c_m:=c(t_m)$ and $(a_i)_2=c_{m+1}:=c(t_{m-1})$.
The index $m$ is chosen so that the signs  of $\delta_e(t_i)-\delta_c$ (for the elastic case) or $\delta_v(t_i)-\delta_c$ (for the visco-elastic case) are different for $a_i=(a_i)_1$ and $a_i=(a_i)_2$.
At the start of crack growth, we begin with $(a_i)_1=a_0=c_0$ and $(a_i)_2=c_1$.
The advantage of choosing previous CZ tip positions, $c_m$, as initial approximations for $a_i$ is that we already know the stress history at these points since they were computed in the previous time steps.

\subsubsection{Computing the Stress at a Crack Tip Position}
Further, during the secant iterations to obtain $a_i$, we will need to compute the stress $\sigma(a_i,t_i)$ for cases when $a_i$ does not equal to the (previous) $c_m$ values. This stress value will be used for the integration while calculating
$K(c_i,t_i)$ and $\delta_e(t_i)$ by \eqref{pes} and \eqref{deltaelas}, respectively.
Note that we cannot directly use the solution given by (\ref{ab}) to solve equation (\ref{pe1}) with $x=a_i$, $c_m<a_i<c_{m+1}$, since it needs $t_c(a_i)$, which approximate calculation can be time consuming.
Instead, we first find $\sigma^{\beta}(c_m,t_i)$ and $\sigma^{\beta}(c_{m+1},t_i)$ by (\ref{ab}) and then employ the following linear interpolant to approximate $\sigma^{\beta}(a_i,t_i)$,
$$
\sigma^{\beta}(a_i,t_i)\approx \sigma^{\beta}(c_m,t_i)+ \frac{\sigma^{\beta}(c_{m+1},t_i)- \sigma^{\beta}(c_m,t_i)}{c_{m+1}-c_m} \left(a_i-c_m\right),\quad
c_m<a_i<c_{m+1}.
$$

\subsubsection{Calculating the Visco-elastic Crack Tip Opening}
The crack tip opening in the visco-elastic case is given by equation (\ref{d=dcvisc2}). However,
when an approximation for $a_i$ is a previous value of $c$, i.e., $a_i=c_m$, then the integration over $\tau$ in equation (\ref{d=dcvisc2}) is from $t_m$ to $t=t_i$.

However, when $a_i\neq c(t_m)$,  the time instant when $a_i$ became part of the CZ, $t_c(a_i)$, is unknown, and to avoid a time consuming calculation of $t_c(a_i)$ we implement the following approach.
%
%
We replace $t_c(a_i)$ with $t_m$, where $c_{m}<a_i<c_{m+1}$, but take into account that $u_e(a_i,t_m)=0$. The integral over $\tau$ will be evaluated by piecewise linearly interpolating $u_e(a_i,\tau)$ between $t_k$ and $t_{k+1}$ for $k=m,m+1,...,i-1$.
Thus, the integral would be written as
$$\sum_{k=m}^{i-1}\int_{t_k}^{t_{k+1}}e^{-\frac{t_i-\tau}{\theta}} [u_e](a_i,\tau)d\tau$$
where
$$[u_e](a_i,\tau)\approx [u_e](a_i,t_k) - \frac{[u_e](a_i,t_{k+1})-[u_e](a_i,t_k)}{t_{k+1}-t_k}(\tau-t_k),\quad
t_k<\tau<t_{k+1},$$
and the first non-zero $[u_e](a_i,t_k)$ is for $k=m+1$, when $c_m<a_i<c_{m+1}$.

During implementation of the algorithm, we come across the step, $i$, where $a_i$ will exceed $c_{i-1}$, and for decreasing CZ length we will have $a_i>c_{i-1}$ in all the steps which follow.
Thus, for these steps, only 1 previous value of $c$ (namely $c_{i-1}$) can be taken as an initial approximation of $a_i$.
To avoid this effect, we will modify the algorithm by fixing $a_i=c_{i-1}$ and computing the corresponding $t_i$ and $c_i$ by solving equation (\ref{d=dc}) (in the elastic case) or (\ref{d=dcvisc}) (in the visco-elastic case) (setting the crack tip opening displacement equal to the critical crack tip opening) and
$K(c_i,t_i)=0$ (setting the stress intensity factor to 0) respectively.

\subsection{Numerical Results for the Propagating Crack}

Choosing PMMA as a reference material, see Section~\ref{RdPMMA},
we used in our numerical examples the value $\delta_c=0.238$  for the normalised critical crack tip opening.

\begin{figure}[H]
\centering
\mbox{\subfigure{\includegraphics[scale=0.55]{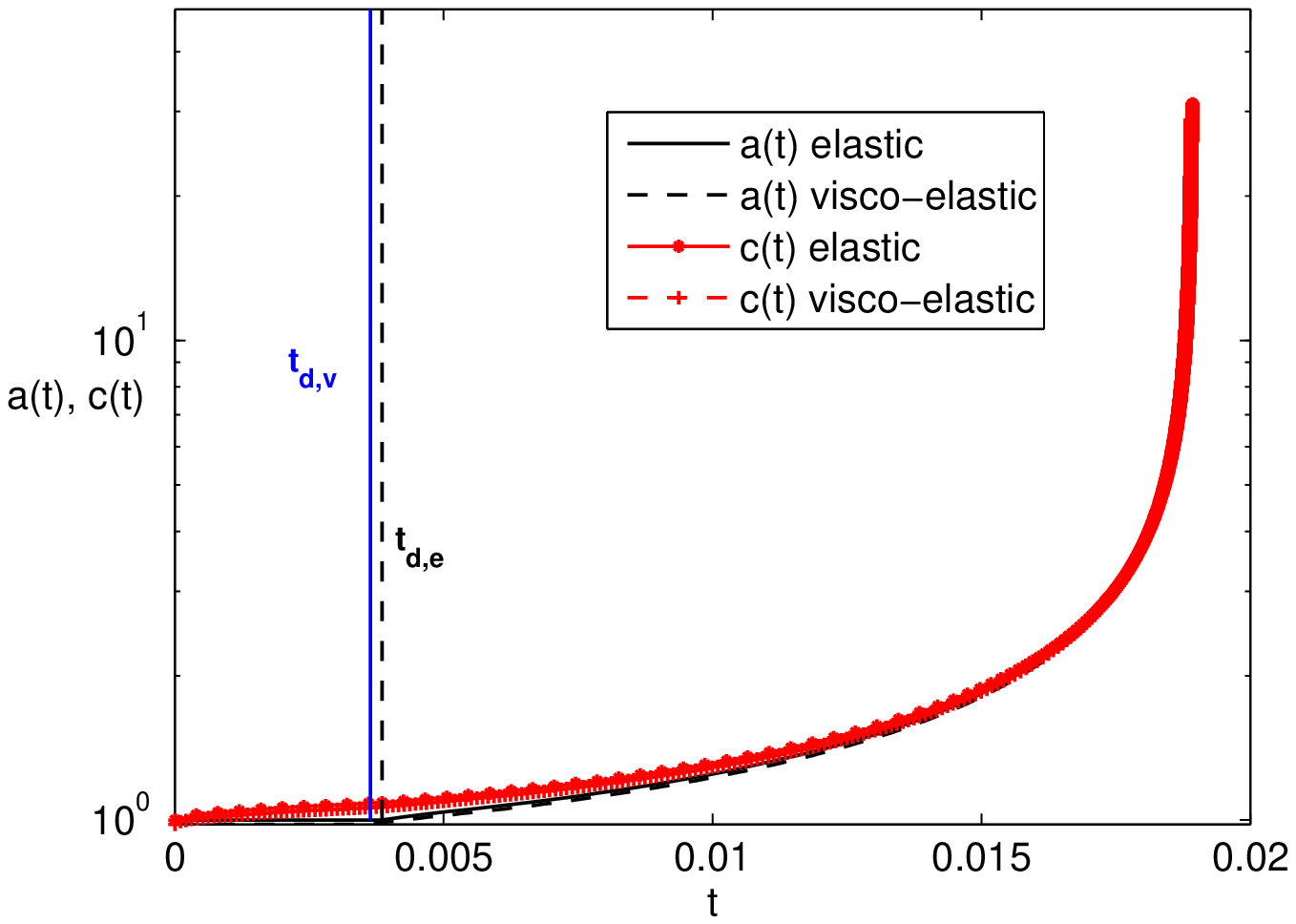}}\quad
\subfigure{\includegraphics[scale=0.55]{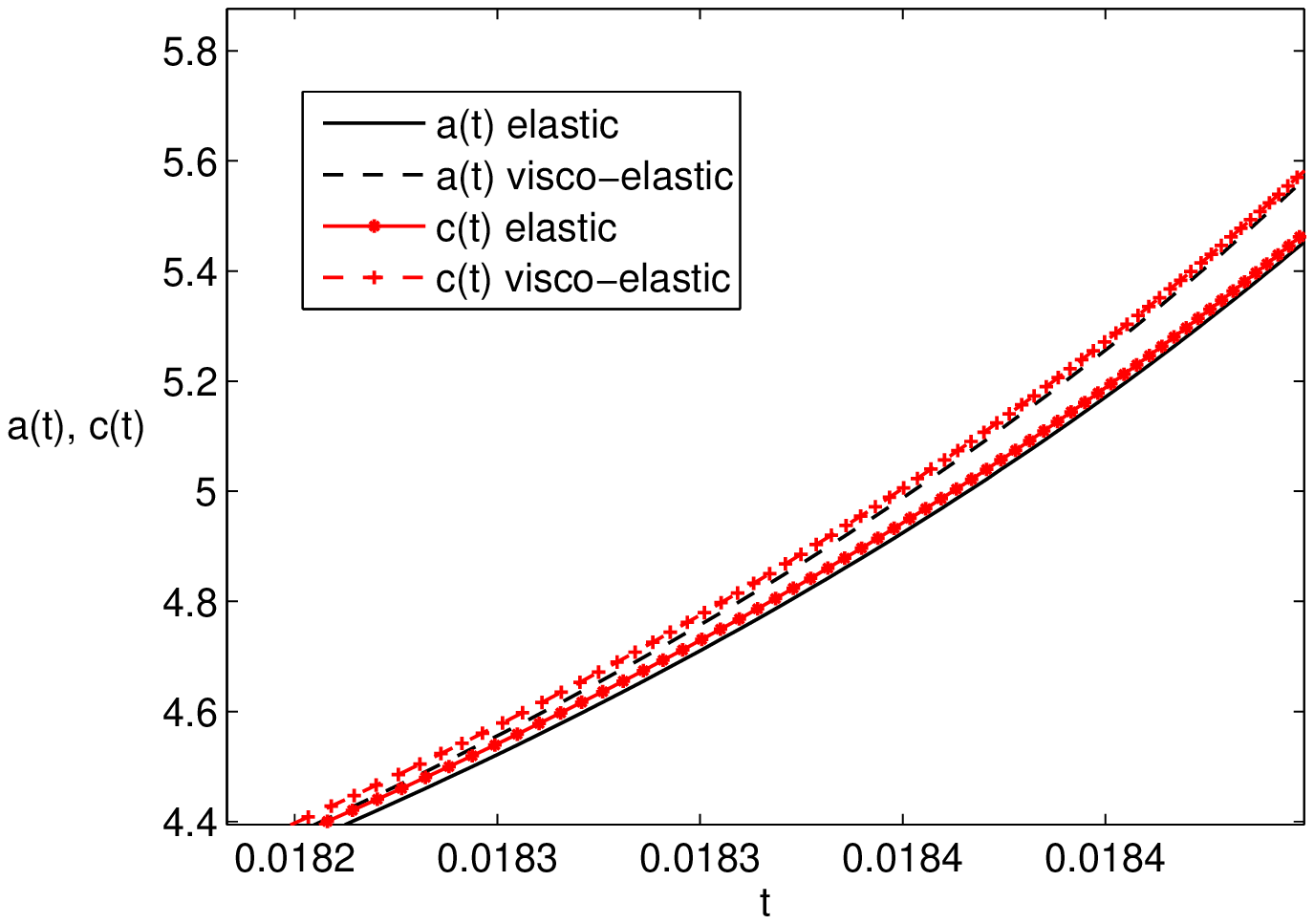}}}
\caption{The crack tip coordinate, $a(t)$, and the CZ tip coordinate, $c(t)$, vs. time, $t$, for $b=4$, $\beta=2$
\\
} \label{cg}
\end{figure}

The graph in Fig. \ref{cg} shows coordinates of the crack tip and the CZ tip for both the elastic and visco-elastic cases.

\begin{figure}[H]
\centering
\includegraphics[scale=0.6]{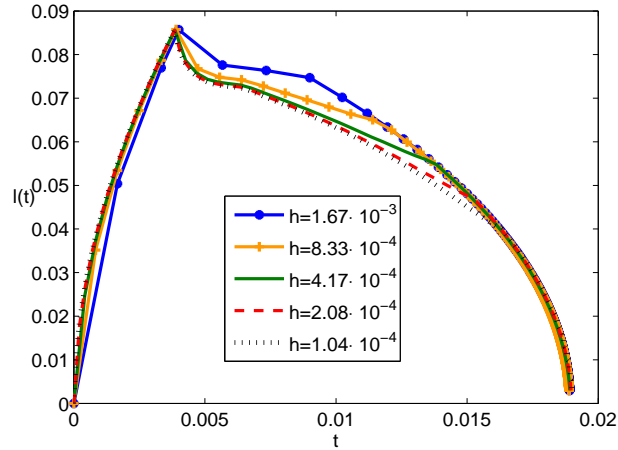}
\captionof{figure}{CZ length, $l(t)$, vs time, $t$, for $b=4$, $\beta=2$ and different time mesh steps, $h$,  (elastic)}
\label{czlb4B2}
\end{figure}

\begin{figure}[H]
\centering
\includegraphics[scale=0.6]{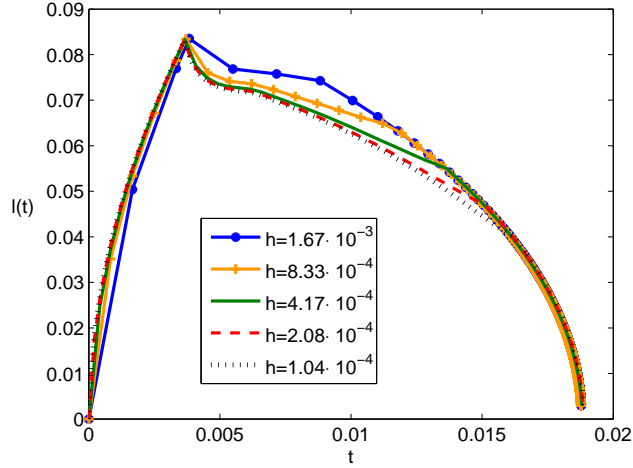}
\captionof{figure}{CZ length, $l(t)$, vs time, $t$, for $b=4$, $\beta=2$ and different time mesh steps, $h$, (visco-elastic)
}
\label{czlb4B2v}
\end{figure}

The graphs in Figs. \ref{czlb4B2} and \ref{czlb4B2v}  show the numerical results for the CZ length evolution in time for the elastic and visco-elastic cases, respectively, calculated with different time meshes.

\begin{figure}[H]
\begin{center}
\includegraphics[scale=0.55]{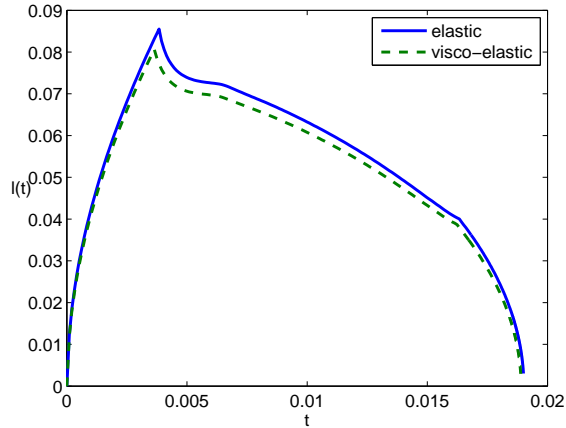}
\caption{CZ length, $l(t)$, vs. time, $t$, for $b=4$, $\beta=2$}\label{cgz2}
\end{center}
\end{figure}

Fig. \ref{cgz2} combines the graphs from Figs. \ref{czlb4B2} and \ref{czlb4B2v} at the finest time-mesh, with $h=4\cdot 10^{-4}$, to compare the evolution of the normalised CZ length, $l(t)=c(t)-a(t)$, with time for the elastic and the visco-elastic cases.  The maxima of these graphs are reached at the delay (crack start) times, $t_d=0.00385$ for the elastic case and $t_d=0.00364$ for the visco-elastic case.

The rate of crack growth, the normalised rupture time, $t_r$, as well as the the CZ length at the crack start time, $l(t_d)$, (which for many cases is the maximum of $l$ in time) depend on the material parameters.
Fig. \ref{mmsrup} and Table~\ref{TT1} show the rupture time dependence of $\beta$ for $b=4$; whereas Fig. \ref{mmsrup2} and Table~\ref{TT2} show the rupture time dependence of $b$ for $\beta=\frac{1}{2}$. The calculations were done for $\beta=\frac{3b}{4}$, $\beta=\frac{b}{2}$, $\beta=\frac{b}{4}$, $\beta=\frac{b}{6}$, and $\beta=\frac{b}{8}$.
The data indicate a strong dependence of the normalised rupture time in the infinite plane on the presence of the crack and on the material parameters $b$ and $\beta$. This is in contrast to the crack propagation results obtained in the models without the cohesive zone, where the rupture time in the plane with and without crack was the same \cite{MikhNamLiv02, MNGeneva03Creep, HakMikIMSE10}.

\begin{minipage}{\textwidth}
\begin{minipage}[c]{0.5\textwidth}
\includegraphics[scale=0.5]{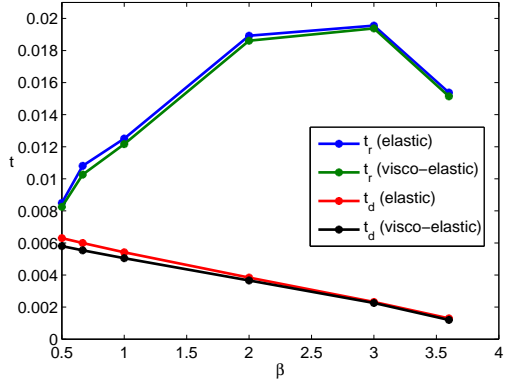}
\captionof{figure}{ Rupture time, $t_r$, vs. $\beta$, for $b=4$ \label{mmsrup}}
  \end{minipage}
  \hfill
  \begin{minipage}[c]{0.4\textwidth}
\begin{tabular}{|c|c|c|}
  \hline
  $\beta$ & $t_{r,e}$ & $t_{r,v}$  \\ \hline
  $3$    & 0.01955  &  0.01938 \\ \hline
  $2$    & 0.01892  &  0.01861 \\ \hline
  $1$    & 0.01251  &  0.01216 \\ \hline
  $2/3$  & 0.01081  &  0.01026 \\ \hline
  $1/2$  & 0.00850  &  0.00825 \\ \hline
\end{tabular}
\captionof{table}{Rupture time for $b=4$, $t_{r,e}$ and $t_{r,v}$ represent the elastic and visco-elastic cases respectively}\label{TT1}
    \end{minipage}
\end{minipage}

\begin{minipage}{\textwidth}
  \begin{minipage}[c]{0.5\textwidth}
    \includegraphics[scale=0.5]{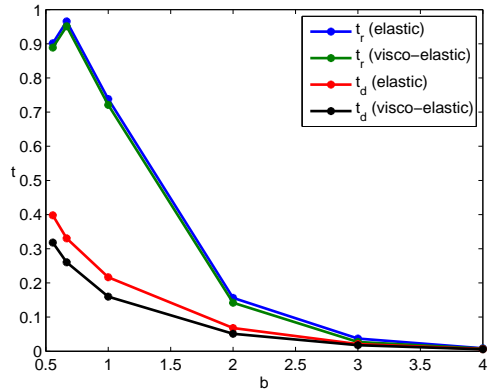}
\captionof{figure}{Rupture time, $t_r$, vs. $b$, for $\beta=1/2$\label{mmsrup2}}
  \end{minipage}
  \hfill
  \begin{minipage}[c]{0.4\textwidth}
\begin{tabular}{|c|c|c|}
  \hline
  $b$ & $t_{r,e}$ & $t_{r,v}$  \\ \hline
  $4$ & 0.00850 & 0.00755 \\ \hline
  $3$ & 0.03714 & 0.0271 \\ \hline
  $2$    & 0.1563 & 0.1418 \\ \hline
  $1$    & 0.7384 &  0.7212 \\ \hline
  $2/3$  & 0.9654 & 0.9513 \\ \hline
\end{tabular}
\captionof{table}{Rupture time for $\beta=1/2$; $t_{r,e}$ and $t_{r,v}$ represent the elastic and visco-elastic cases respectively}\label{TT2}
\end{minipage}
\end{minipage}

\ \\
The dependency of the CZ length at $t=t_d$ for different parameter sets, is given in Figs. \ref{mmstdnew} and \ref{mmstdnew2} as well as in Tables \ref{TT1l} and \ref{TT2l}. As one can see from Fig.~\ref{mmstdnew2}, the CZ length can reach maximum not at $t=t_d$ but after the crack start, for small $b$ in the elastic case.

\begin{minipage}{\textwidth}
\begin{minipage}[c]{0.5\textwidth}
\includegraphics[scale=0.5]{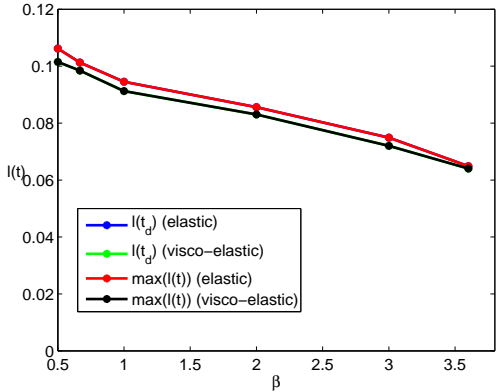}
\captionof{figure}{CZ length at $t=t_d$, vs. $\beta$, for $b=4$ \label{mmstdnew}}
  \end{minipage}
  \hfill
  \begin{minipage}[c]{0.4\textwidth}
\begin{tabular}{|c|c|c|}
  \hline
  $\beta$ & $l(t_{d,e})$ & $l(t_{d,v})$  \\ \hline
  $3$    & 0.07490  &  0.07203 \\ \hline
  $2$    & 0.08562  &  0.08304 \\ \hline
  $1$    & 0.09447  &  0.09126 \\ \hline
  $2/3$  & 0.1013  &  0.09843 \\ \hline
  $1/2$  & 0.1062  &  0.10150 \\ \hline
\end{tabular}
\captionof{table}{Elastic and visco-elastic CZ length at $t_{d,e}$ and $t_{d,v}$, respectively, vs. $\beta$, for $b=4$.}\label{TT1l}
    \end{minipage}
\end{minipage}

\begin{minipage}{\textwidth}
\begin{minipage}[c]{0.5\textwidth}
\includegraphics[scale=0.5]{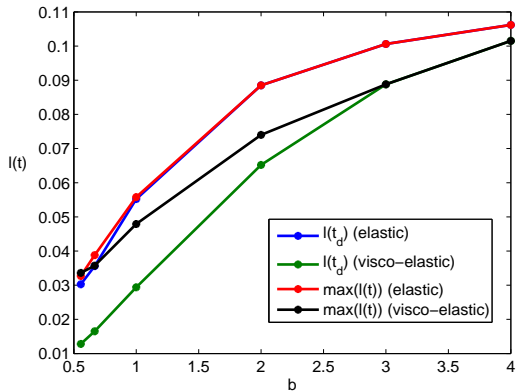}
\captionof{figure}{CZ length at $t=t_d$ and maximal CZ length, vs. $b$, for $\beta=\frac{1}{2}$ \label{mmstdnew2}}
  \end{minipage}
  \hfill
  \begin{minipage}[c]{0.4\textwidth}
\begin{tabular}{|c|c|c|}
  \hline
  $b$ & $l(t_{d,e})$ & $l(t_{d,v})$  \\ \hline
  $4$    & 0.1062  &  0.1015 \\ \hline
  $3$    & 0.1006  &  0.0888 \\ \hline
  $2$    & 0.0885  &  0.0652 \\ \hline
  $1$    & 0.0552  &  0.0294 \\ \hline
  $2/3$  & 0.0356  &  0.0165 \\ \hline
\end{tabular}
\captionof{table}{CZ length at the delay time for $\beta=\frac{1}{2}$, $t_{d,e}$ and $t_{d,v}$ represent the delay times of the elastic and visco-elastic cases respectively}\label{TT2l}
    \end{minipage}
\end{minipage}

\subsubsection{Crack Start Jump Analysis}
\begin{figure}[H]
\begin{center}
\includegraphics[scale=0.55]{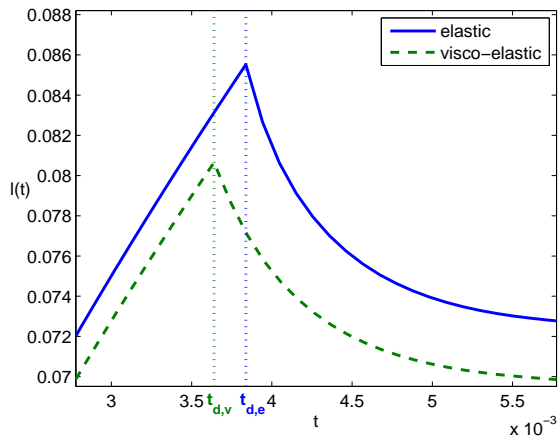}
\caption{CZ length, $l(t)$, vs. time, $t$, for $b=4$, $\beta=2$ \label{cgm}}
\end{center}
\end{figure}

Let us now analyse what is happening at the onset on crack growth.
Considering the graph of the CZ length, $l(t)$, at the vicinity of the delay time, see Fig. \ref{cgm},
we see that for $b=4$ and $\beta=1/2$, the CZ length $l(t)$ reaches maximum at the corresponding values of $t_d$ (different for the elastic and visco-elastic cases),
and the decrease of the CZ length as crack growth begins is more pronounced in the elastic case compared with the visco-elastic one.
Moreover,
by taking other sets of parameters, $b=4$, $\beta=1/2$, we observe that  the decrease of $l(t)$ at the onset of crack growth is sharper.

\begin{figure}[H]
\centering
\mbox{\subfigure{\includegraphics[scale=0.55]{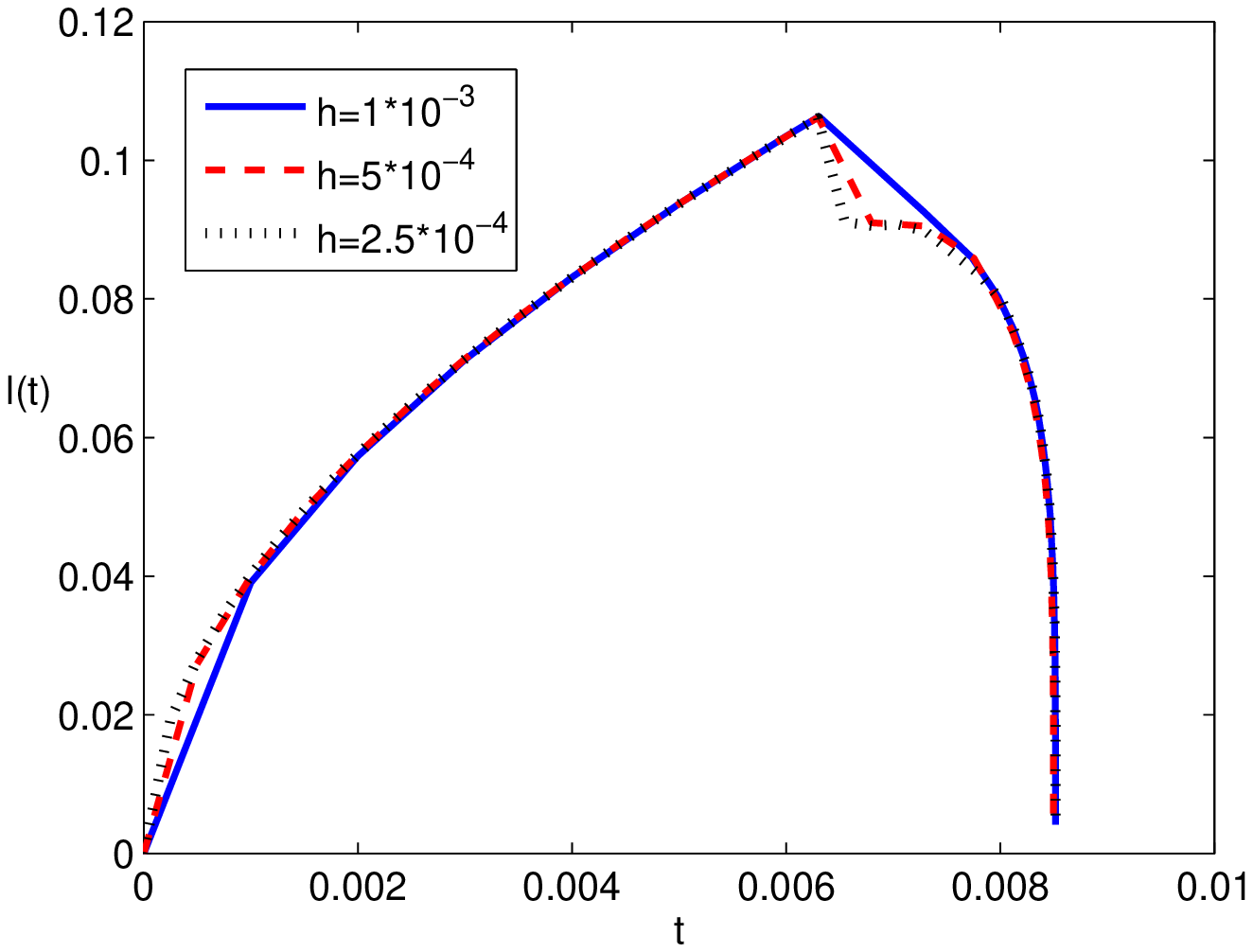}}\quad
\subfigure{\includegraphics[scale=0.55]{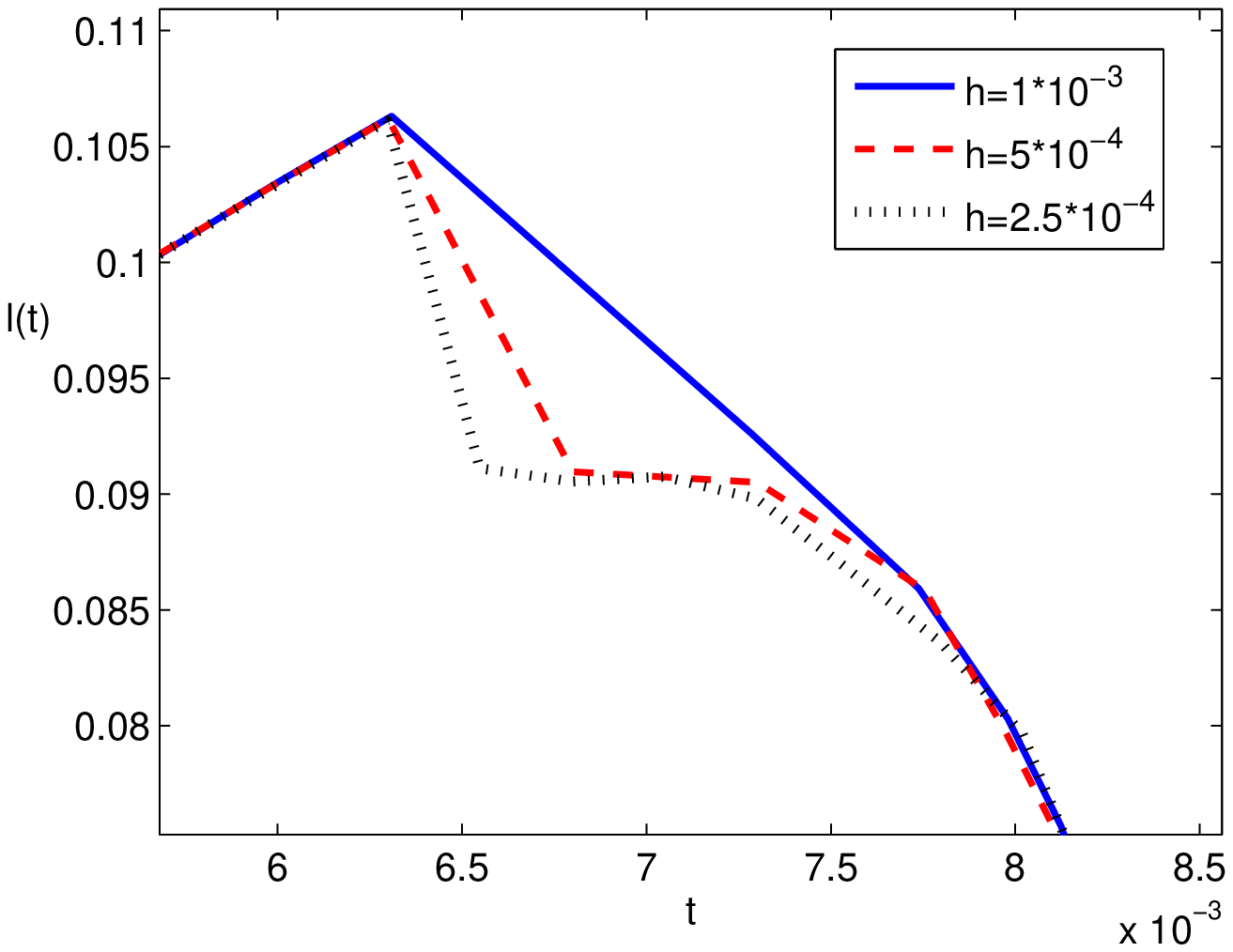}}}
\caption{CZ length, $l(t)$, vs time, $t$, for  $b=4$, $\beta=1/2$ (elastic)}\label{fb11m}
\end{figure}

Let us consider, how changing the time step size influences the CZ length behaviour obtained numerically.
We have 7, 13, and 26 time steps before crack growth for cases $h=1\cdot 10^{-3}$, $h=5\cdot 10^{-4}$, and $h=2.5\cdot 10^{-4}$ respectively. From Fig.~\ref{fb11m} we see that the initial decrease in the CZ length is sharper at finer time meshes (and hence more steps before the crack growth start).

\begin{figure}[H]
\centering
\includegraphics[scale=0.55]{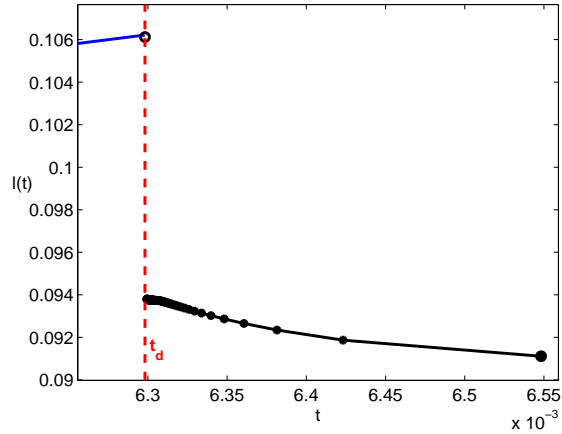}
\caption{CZ length, $l$, vs. time, $t$, near $t_d$ for $b=4$, $\beta=1/2$ (elastic)
\\
}
\label{mmstria}
\end{figure}

\begin{figure}[H]
\begin{center}
\includegraphics[scale=0.55]{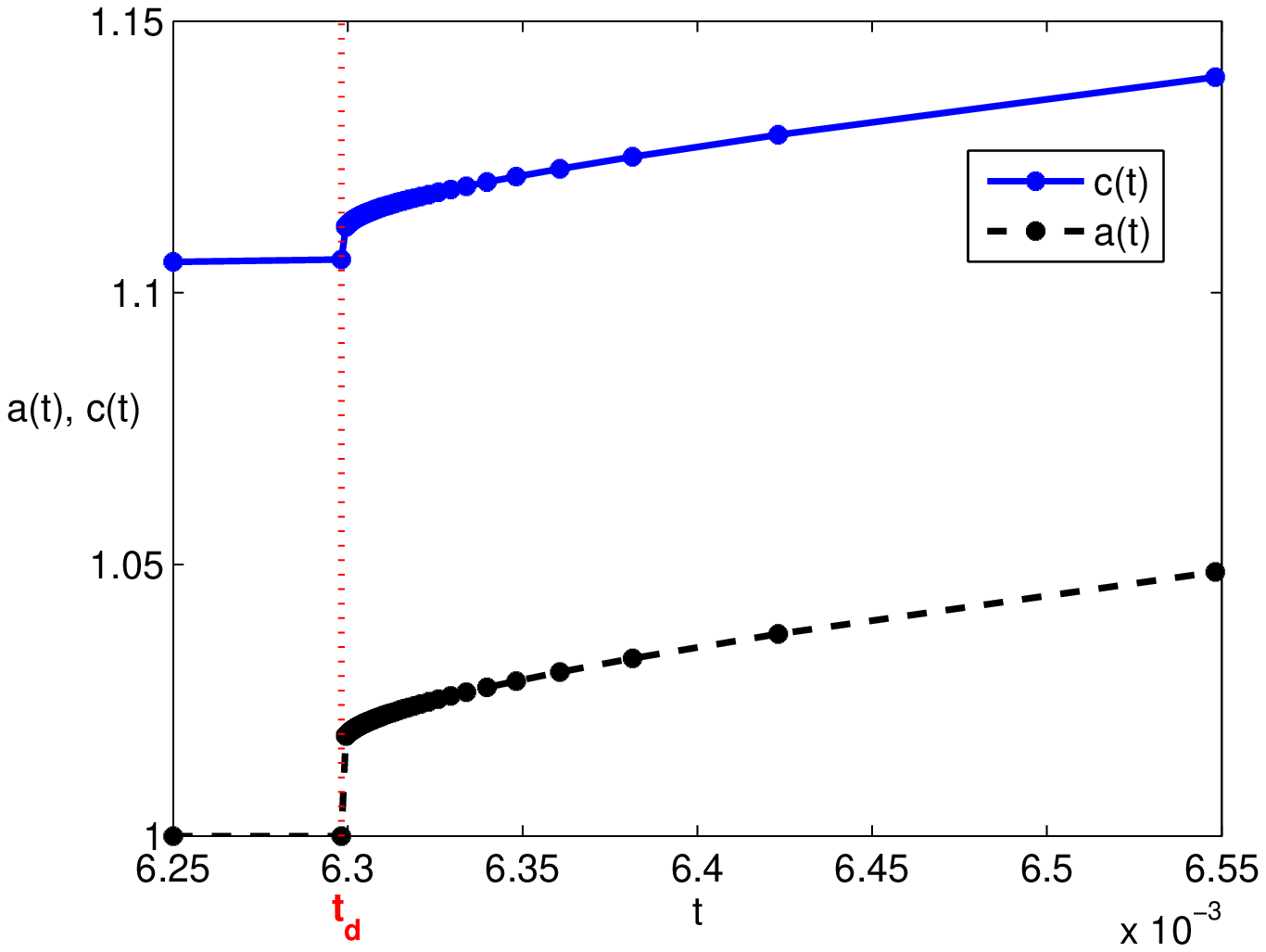}
\caption{$a(t)$ and $c(t)$  vs. time, $t$, for $b=4$, $\beta=1/2$ (elastic)}\label{mmsextria2}
\end{center}
\end{figure}

\begin{figure}[H]
\begin{center}
\includegraphics[scale=0.55]{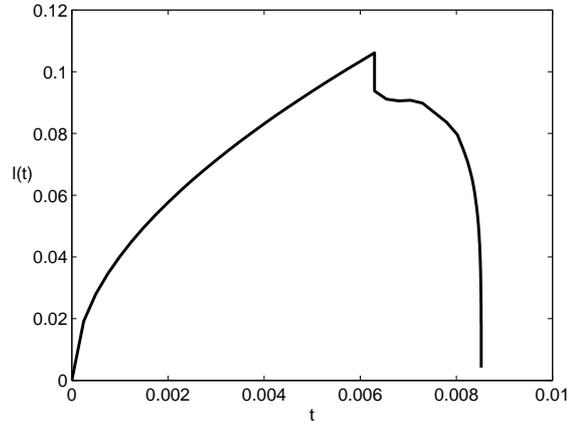}
\caption{CZ length, $l(t)$ vs. time, $t$, for $b=4$, $\beta=1/2$}\label{vj}
\end{center}
\end{figure}

To analyse, in a more systematic way, whether the CZ tip coordinate, $c(t)$, crack tip coordinate, $a(t)$, and the CZ length, $l(t)=c(t)-a(t)$, are continuous or discontinuous at $t=t_d$,  we calculated there vales at a sequence of points $t_{i}>t_d$ tending to $t_d$. The numerical experiment for $b=4$, $\beta=1/2$, see Figs. \ref{mmstria}, \ref{mmsextria2}, and \ref{vj}, at $h=2.5\cdot 10^{-4}$,
shows that $a(t_i)$, $c(t_i)$ and $l(t_i)$ tend, respectively to the limiting values $a_{d+}= 1.018\not=a_0=1$,  $c_{d+}= 1.112\not=c(t_d)=1.106$, and $l_{d+}=0.0938\not=l(t_d)= 0.106$.

{\em The jump of the crack length at $t=t_d$, seen on the figures, indicates that there is an unstable crack growth at the onset of crack propagation followed by the stable crack growth,  for the chosen set of parameters, $b=4$, $\beta=1/2$. It causes also a jump decrease in the CZ length followed by  a continuous CZ length evolution.}

\begin{figure}[H]
\begin{center}
\includegraphics[scale=0.55]{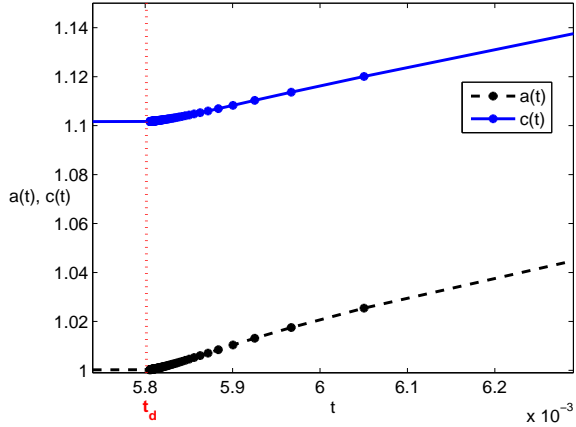}
\caption{$a(t)$ and $c(t)$ vs. time, $t$, for $b=4$, $\beta=1/2$ (visco-elastic)}\label{mmstriav}
\end{center}
\end{figure}

A similar analysis for the visco-elastic case with the same parameters, $b=4$, $\beta=1/2$, see Fig.~\ref{mmstriav}, shows that the functions
$c(t)$, $a(t)$, and $l(t)$ are continuous at $t=t_d$, unlike in the elastic case.

\section{Convergence Rates}\label{Conv}

In the model problems which were numerically solved, we obtained numerical solutions using successively refined meshes.
Now we will look at the convergence rate of several computed variables in more details.

Let $y$ denote the exact  value of a variable, and
$y_N$ corresponds to the numerical solution obtained for the step size $h=h_N$, $N=1,2,...,{N_*},$  and
$$\epsilon=\epsilon(y_N)=|y-y_N|$$
denote the corresponding absolute error.

When the exact value are unknown,
we will use Aitken's extrapolation technique, also known as the Aitken $\Delta^2$ process,
see for example \cite[Section 2.6]{{Atkinson1998}}, to accelerate the
convergence and obtain a good approximation to the exact solution.
It gives an approximation $y_a$ of the exact solution $y$ based on 3
consecutive terms of a convergent sequence,
\be
y_a=
\frac{y_{_{N_*}}y_{_{{N_*}-2}}-\left(y_{_{{N_*}-1}}\right)^2}{y_{_{N_*}}-2y_{_{{N_*}-1}}+y_{_{{N_*}-2}}}
\label{AEF}.
\ee

Consequently, the approximate error will be taken as
$$\epsilon_N\approx|y_a-y_N|.$$

We assume that
there exists a constant $C$ such that
$\epsilon_N=Ch_N^{\alpha}$.
Then we have
$$\alpha=\frac{\log\left(\epsilon_{N-1}/\epsilon_{N}\right)}{\log\left(h_{N-1}/h_{N}\right)}.$$

\subsection{Convergence for Stationary Crack}
Taking the time instant $t=0.6$,
we plot in Fig.~\ref{mmscratef1} the error in the CZ length, $l$, for $b=4$, $\beta=2$, verses the time step, $h$, and present the numerically estimated convergence rate in Table~\ref{T1}.

\begin{minipage}{\textwidth}
  \begin{minipage}[c]{0.6\textwidth}
    \centering
    \includegraphics[scale=0.5]{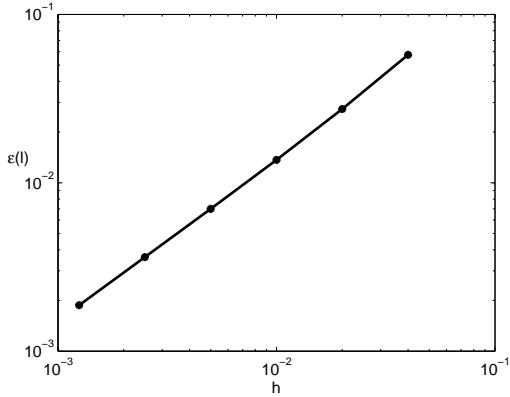}
    \captionof{figure}{CZ length error, $\epsilon(l)$, vs time step, $h$, \\for $b=4$, $\beta=b/2$, $t=0.6$\\}\label{mmscratef1}
  \end{minipage}
  \hfill
  \begin{minipage}[c]{0.4\textwidth}
    \centering
\begin{tabular}{|c|c|}
  \hline
  $h$ & $\alpha$  \\ \hline
  $0.04$& -  \\ \hline
  $0.02$&  1.06782\\ \hline
  $0.01$ & 1.00550 \\ \hline
  $0.005$ &0.96689 \\ \hline
  $0.0025$&0.94986 \\ \hline
  $0.00125$&0.94986\\ \hline
\end{tabular}
      \captionof{table}{\ \\ Convergence rate for CZ length, $l$ }\label{T1}
    \end{minipage}
  \end{minipage}

We are particularly interested in convergence of stress at the CZ tip. The following table shows the values of $\sigma$ at the CZ tip in the elastic case at $t=0.6$ for $b=4$, $\beta=2$ and $\beta=1/2$, obtained for different time steps, $h$.
\begin{table}[H]
\begin{center}
\begin{tabular}{|c|l|l|l|l|l|l|}
  \hline
  $h$&0.04 & 0.02 & 0.01 & 0.005 & 0.0025 & 0.00125 \\ \hline
  $\sigma|_{\beta=2}$&1.39334 & 1.44340 &  1.49319 & 1.54250 & 1.59121 &  1.63929  \\\hline
  $\sigma|_{\beta=1/2}$&1.17767 & 1.18144 &  1.18438 & 1.18667 & 1.18844 & 1.18980 \\
  \hline
\end{tabular}
\end{center}
\end{table} %
Using Aitken's extrapolation formula \eqref{AEF}, we obtain $\sigma_a=5.25425$ at $\beta=2$, and $\sigma_a=1.19443$ at $\beta=1/2$, as the approximations to the exact solutions.
The error, $\epsilon(\sigma)$, is presented in Fig.~\ref{mmscrates1}, and the corresponding approximate convergence rates, $\alpha|_{\beta=2}$ and $\alpha|_{\beta=1/2}$, are given in Table~\ref{T2}.

\begin{minipage}{\textwidth}
  \begin{minipage}[c]{0.45\textwidth}
    \includegraphics[scale=0.5]{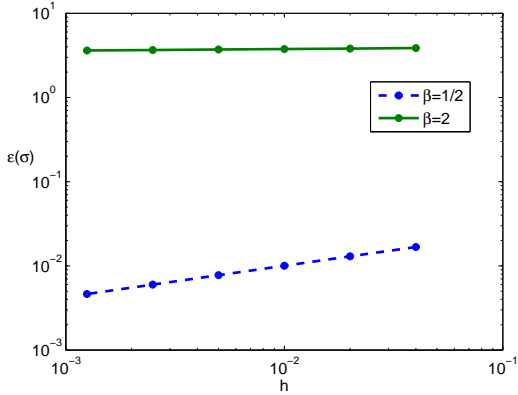}
\captionof{figure}{\ \\ CZ tip stress error, $\epsilon(\sigma)$,  vs time step, $h$,\\
 for $b=4$;{\tiny } $\beta=2$ and $\beta=1/2$, at $t=0.6$\\}
 \label{mmscrates1}
  \end{minipage}
  \hfill
  \begin{minipage}[c]{0.4\textwidth}
\begin{tabular}{|c|c|c|}
  \hline
  $h$ & $\alpha|_{\beta=2}$ & $\alpha|_{\beta=1/2}$  \\ \hline
  $0.04$    & -  & - \\ \hline
  $0.02$    & 0.01883 & 0.3675 \\ \hline
  $0.01$    & 0.01897 & 0.3708 \\ \hline
  $0.005$   & 0.01904 & 0.3724 \\ \hline
  $0.0025$  & 0.01906 & 0.3730 \\ \hline
  $0.00125$ & 0.01906 & 0.3730 \\ \hline
\end{tabular}
      \captionof{table}{\ \\Convergence rates for the CZ tip stress,\\
      $\sigma$,
      for $\beta=2$ and $\beta=1/2$}\label{T2}
    \end{minipage}
\end{minipage}

The very slow (if at all) convergence of the stress at the CZ tip, when $b=4$ and $\beta=2$, cf. also Fig.~\ref{fig:test2}, may be a manifestation of a CZ tip stress singularity in the exact solution, at some range of parameters $b$ and $\beta$, in spite we assumed that the stress is bounded at the CZ tip. An a priory information about the stress singularity in the considered nonlinear problem would be useful but is not available. Constructing of stress asymptotics is beyond the scope of this paper but may be considered elsewhere.

Fig.~\ref{crateu1} and Table~\ref{T5} show  the error, $\epsilon(\delta)$, and the numerical convergence rates, $\alpha_e$ and $\alpha_v$, of the crack tip opening, $\delta$, for the elastic and the visco-elastic cases, respectively, versus the time step, $h$.

\noindent
\begin{minipage}{\textwidth}
  \begin{minipage}[c]{0.6\textwidth}
    \includegraphics[scale=0.5]{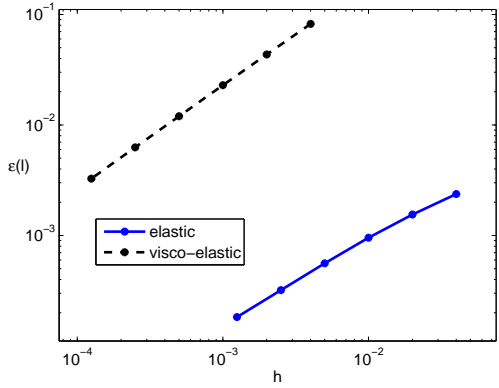}
\captionof{figure}{Crack tip opening error, $\epsilon(\delta)$,\\
vs time step, $h$,
for $b=4$, $\beta=2$
}\label{crateu1}
  \end{minipage}
  \hfill
  \begin{minipage}[c]{0.4\textwidth}
\begin{tabular}{|c|c|c|}
  \hline
  $h$ & $\alpha_e$ & $\alpha_v$  \\ \hline
  $0.04$& -    & -  \\ \hline
  $0.02$& 0.6176 & 0.9176 \\ \hline
  $0.01$& 0.6963 & 0.9234 \\ \hline
  $0.005$&0.7692 & 0.9321 \\ \hline
  $0.0025$& 0.8087&0.9381  \\ \hline
  $0.00125$& 0.8087&0.9381 \\ \hline
\end{tabular}
      \captionof{table}{\ \\Convergence rates for crack tip\\ opening,  $\delta$}\label{T5}
    \end{minipage}
\end{minipage}

\subsection{Convergence for Propagating Crack}


Now, we will compute the convergence rates at a time instants before and after the crack start time $t_d$.
Note that the CZ length evolution is the same for the elastic and visco-elastic cases if $t<t_d$. We show the results for $b=4$, $\beta=2$, for which, as was mentioned above, $t_d=0.00385$ for the elastic case and $t_d=0.00366$ for the visco-elastic case.

Let $\alpha_1$ denote the numerical convergence rate at $t=1/600\approx 0.0017<t_d$, while $\alpha_{2e}$ and $\alpha_{2v}$ denote the numerical convergence rates, for the elastic and visco-elastic cases respectively, at $t=0.01>t_d$. Fig.~\ref{mmscrateczcomp} presents the graphs of the errors, while Table~\ref{T6} shows the order of convergence, of the CZ length.

\begin{minipage}{\textwidth}
  \begin{minipage}[c]{0.5\textwidth}
    \includegraphics[scale=0.5]{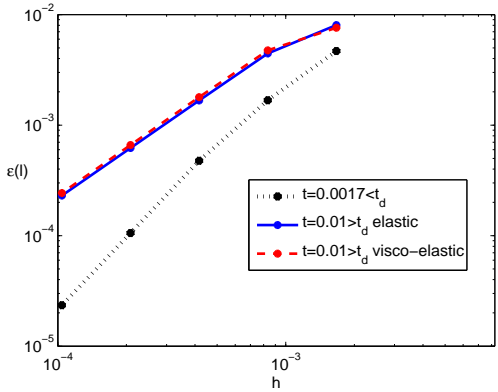}
\captionof{figure}{\ \\CZ length error, $\epsilon(l)$, vs time step, $h$, for $b=4$, $\beta=2$
}\label{mmscrateczcomp}
  \end{minipage}
  \hfill
  \begin{minipage}[c]{0.5\textwidth}
\begin{tabular}{|c|c|c|c|}
  \hline
  $h$ & $\alpha_1$ & $\alpha_{2e}$ & $\alpha_{2v}$  \\ \hline
  $0.0017$& -& -&-  \\ \hline
  $0.00083$& 1.47845&  0.85015& 0.67726\\ \hline
  $0.00042$ & 1.82147& 1.41520& 1.40615\\ \hline
  $0.00021$ & 2.16984& 1.42908& 1.44186\\ \hline
  $0.00010$& 2.16984&  1.42908& 1.44186\\ \hline
\end{tabular}
\captionof{table}{\ \\Convergence rates for the CZ length, $l$}\label{T6}
    \end{minipage}
\end{minipage}

\section{Some remarks, auxiliary proofs, and material parameters}
\subsection{Material parameter range of CZ model applicability}\label{AP-CZ}
To analyse the range of $\gamma=\beta/b$, for which the CZ model based on conditions \eqref{czc1}, \eqref{czc} can exist, let us first remark that if $\gamma=1$, then
\be
\frac{d}{d\hat t}\underline{\Lambda}(\boldsymbol{\hat{\sigma}};\hat{t})
=\frac{1}{b \sigma_0^{\beta}}\underline{\Lambda}^{1-\beta}(\boldsymbol{\hat{\sigma}};\hat{t}) |\boldsymbol{\hat{\sigma}}(\hat{t})|^\beta,
\nonumber\ee
which means that $\underline{\Lambda}(\boldsymbol{\hat{\sigma}};\hat{t})$ is a strictly growing function for any $\hat t$, when $|\boldsymbol{\hat{\sigma}}(\hat{t})|>0$. Particularly, if the cohesive condition \eqref{czc} is reached at some point $\hat x$ in time $\hat t_c(\hat x)$, it can not stay at larger times, $\hat t>\hat t_c(\hat x)$ unless $\boldsymbol{\hat{\sigma}}(\hat x,\hat{t})=0$; but if $|\boldsymbol{\hat{\sigma}}(\hat x,\hat{t})|=0$ for $t>t_c(\hat x)$ this means that the point $\hat{x}$ belongs to the crack rather than to the CZ. That is, the CZ can not exist at $\gamma=1$ (which corresponds to the Robinson damage linear accumulation rule). Instead, the Robinson damage linear accumulation rule implies the crack propagation without the cohesive zone, which, as follows from \cite{MikhNamLiv02, MNGeneva03Creep, HakMikIMSE10}, is possible only if $0<b<2$ in the local approach, although the Neuber-Novozhilov type non-local approach extends the applicability range to arbitrary $b>0$.

Similarly, if $\gamma>1$, then
\be
\frac{d}{d\hat t}\underline{\Lambda}(\boldsymbol{\hat{\sigma}};\hat{t})
=\frac{\gamma}{\beta \sigma_0^{\beta}}\underline{\Lambda}^{1-\beta}(\boldsymbol{\hat{\sigma}};\hat{t}) \int_0^{\hat{t}} |\boldsymbol{\hat{\sigma}}(\hat{\tau})|^\beta(\hat{t}-\hat{\tau})^{\gamma-2}d\hat{\tau}>0,
\nonumber\ee
which also prevents for the CZ condition \eqref{czc} to hold at any time $\hat t>\hat t_c(\hat x)$, after the condition had been reached at a time $\hat t_c(\hat x)$, even if $\sigma(\hat x,\hat t)=0$ for $\hat t>\hat t_c(\hat x)$.

Thus, the CZ model defined by \eqref{czc1}, \eqref{czc} is applicable only if $0<\gamma<1$, i.e., $0<\beta<b$. Note that the material parameters obtained in \cite{MikNam2011AAM} by fitting experimental data for several structural materials, satisfy these conditions.

\subsection{Analytical Solution of the Abel Type Equation in the Cohesive Zone}\label{Sab}
To obtain the stresses in the CZ, we have to solve the Abel-type linear integral equation \eqref{pe1} for $\sigma^{\beta}(x,t)$ at $t\ge t_c$, when $\sigma^{\beta}(x,\tau)$ is known at $\tau\in [0,t_c(x)]$  at its right hand side.

To this end, we have the following important assertion, see for example \cite[Theorem 1.2.1]{abelbook}.
\begin{thm}
If $f(t)$ is absolutely continuous on $[t_c,T_1]$, then the Abel type integral equation
\be\int_{t_c}^t g( \tau)(t-\tau)^{\gamma-1}d\tau=f(t),\quad t\in [t_c,T_1],\quad \gamma\in(0,1)\nonumber\ee
has a unique solution  $g$ in $L_1(t_c,T_1)$, which is given by formula
\be g( \tau) = \frac{\sin{(\pi\gamma)}}{\pi} \frac{d}{d\tau}\int_{t_c}^\tau {f(t)}{(\tau-t)^{-\gamma}}dt\label{sol1}.\ee
\end{thm}
Integrating by parts, expression (\ref{sol1}) can be written as
\be g( \tau) = \frac{\sin{(\pi\gamma)}}{\pi}\left(f({t_c})(\tau-{t_c})^{-\gamma} +\int_{{t_c}}^\tau f'(t)(\tau-t)^{-\gamma}dt\right).\label{sol2}\ee

For equation (\ref{pe1}),
$$f(t)=\frac{1}{\gamma} -\int_{0}^{t_c} \sigma^{\beta}(x,\tau)(t-\tau)^{\gamma -1}d\tau.$$
Moreover, we know that $f(t_c)=0$ when $t_c>0$, since the condition $\Lambda=1$ (see equation (\ref{pe1})) is satisfied at $t=t_c>0$. Finally, $\sigma^\beta(x,\tau)=g( \tau)$ at $\tau\ge t_c$.

\subsection{Continuity of $\sigma(x,t)$ in $t$.}
Let us now analyse the behaviour of the numerical solution of (\ref{pe1}) for $\sigma(x,t)$ as $t\to t_c(x)+0$ and prove that if we take the piece-wise approximation of the function $\sigma^\beta(x,t)$ in $t\le t_c(x)$ over the time instants $t_j$, $j=0,1,2,...,k$, we obtain continuity of $\sigma^\beta(x,t)$ at $t=t_c(x)$ also from the right, i.e., when $t\to t_c(x)+0$.

Indeed, from the first equality in \eqref{ab}, where $t_k=t_c(x)$, we have,
\begin{align}&\lim_{t\to t_c(x)+0}\sigma^{\beta}(x,t) = \lim_{t\to t_c(x)+0} \frac{-1}{\pi}\sin{\left(\pi\gamma \right)} \left[ \sum_{j=1}^k\sigma^{\beta}(x,t_{j-1})\left(V(t_{j-1},t,t_c(x))-V(t_j,t,t_c(x))\right)
\right.\nonumber\\ &\left.
+ \frac{1}{\gamma}\left(\frac{\sigma^{\beta}(x,t_j)- \sigma^{\beta}(x,t_{j-1})}{t_j-t_{j-1}}\right)\left(W(t_{j-1},t,t_c(x))-W(t_j,t,t_c(x))- \gamma (t_j-t_{j-1})V(t_j,t,t_c(x))\right) \right].\label{lims}\end{align}
In the above formula, only $V(y,t,t_c(x))$ and $W(y,t,t_c(x))$ depend on $t$. Moreover, since $\beta>0$, from \eqref{Ve} and \eqref{We} we have for any $t$,
$$
V(t_c,t,t_c)=\pi\csc{\left(\pi\gamma \right)},\quad W(t_c,t,t_c)=0.
$$

For the case when $y\not= t_c$, we have
\begin{align}\lim_{t\to t_c+0}V(y,t,t_c)&=\lim_{t\to t_c+0} \left\{\pi\csc{\left(\pi\gamma \right)} -\frac{1}{\gamma}\left(\frac{t_c-y}{t-y}\right)^{\gamma }\,\, _2F_1 \left[\gamma ,\gamma ;1+\gamma ;\frac{t_c-y}{t-y}\right]\right\}
\nonumber\\
&= \pi\csc{\left(\pi\gamma \right)} -\frac{1}{\gamma}\Gamma\left[1+\gamma \right]\Gamma\left[1-\gamma \right]
\nonumber\\ &= \pi\csc{\left(\pi\gamma \right)}-\frac{1}{\gamma}\left(\pi\gamma \csc{\left(\pi\gamma \right)}\right) = 0, \nonumber
\end{align}
\begin{align}&\lim_{t\to t_c+0}W(y,t,t_c)\nonumber\\&=\lim_{t\to t_c+0} \left\{\gamma \pi\csc{\left(\pi\gamma \right)}(t-y)-\frac{1}{1+\gamma}(t_c-y)^{1+\gamma }(t-y)^{-\gamma } \,\,_2F_1 \left[1+\gamma ,\gamma ;2+\gamma ;\frac{t_c-y}{t-y}\right]\right\}
\nonumber\\ & = \gamma\pi \csc{\left(\pi\gamma \right)}(t_c-y)-\frac{1}{1+\gamma}(t_c-y)\Gamma\left[2+\gamma \right]\Gamma\left[1-\gamma \right]
\nonumber\\ & = \gamma\pi \csc{\left(\pi\gamma \right)}(t_c-y)-\gamma\pi(t_c-y)\csc{\left(\pi\gamma \right)} \nonumber = 0,
\end{align}
where we have used that $$_2F_1[a,b,c,1]=\frac{\Gamma[c]\Gamma[c-a-b]}{\Gamma[c-a]\Gamma[c-b]},$$ as well as other properties of the Gamma function such as $$\Gamma[z+1]=z \Gamma[z] \quad \textrm{and} \quad \Gamma[1-z]\Gamma[z]=\pi\csc[(\pi z)].$$

Consequently, in equation (\ref{lims}), the summation over $j$   yields
\begin{multline}\lim_{t\to t_c(x)+0}\sigma^{\beta}(x,t)= -\frac{1}{\pi}\sin{\left(\pi\gamma \right)} \left[ -\pi\csc{\left(\pi\gamma \right)}\sigma^{\beta}(x,t_{k-1}) +
\right.\\
 \left. \frac{1}{\gamma}\cdot\frac{\sigma^{\beta}(x,t_c(x))- \sigma^{\beta}(x,t_{k-1})}{t_c(x)-t_{k-1}} \left(\gamma (t_c(x)-t_{k-1})\pi\csc{\left(\pi\gamma \right)} \right)\right] =\sigma^{\beta}(x,t_c(x)).\nonumber\end{multline}
Therefore, $\lim_{t\to t_c(x)+0}\sigma^{\beta}(x,t)=\sigma^{\beta}(x,t_c(x))$ for  $0<\gamma<1$.

\subsection{Material parameters used in the numerical examples}\label{RdPMMA}
We used in this paper the material parameters close to the ones for PMMA.

For the rheological parameters we took (cf. \cite[pages 655-657]{pdatah}, \cite{viscos}, and \cite{dc}): Poisson's ratio $\nu=0.35$;
Young's modulus of elasticity  $E_0=3100$MPa (hence $\mu_0=1148$MPa);
viscosity $\eta=2 \cdot 10^{7}$ MPa s. We also chosen $\hat\theta = 3.23\cdot 10^4$ s.

Fitting the static creep rupture data under tensile stress for PMMA from \cite{mck},  gives the values $b=18.5$ and $\sigma_0=58.1$MPa\ hr$^{1/b}$ in the durability curve \eqref{dc}. Taking the applied load $\hat{q}=51.6$MPa, we arrive at the values $\hat t_\infty=8.96$hr, $m=5$ and $\theta=1$.

We took the critical crack opening displacement $\hat{\delta}_c=0.0016$mm, cf. \cite[Section 10.3.2]{dc} and references therein.
Under the plane stress condition $\varkappa=(3-\nu)/(1+\nu)=1.96$ and by \eqref{dc} we obtain $\delta_c=0.238$ for $\hat{q}=51.6$MPa and for the initial crack length $\hat{a}_0=0.1$mm.

\section{Concluding remarks}

A novel non-linear history-dependent cohesive zone model of crack propagation in linearly elastic and visco-elastic materials, which is a history-dependent modification of the Leonov-Panasyuk-Dugdale model, was introduced in the paper.
The normal stress on the cohesive zone satisfies the history dependent yield condition, given in terms of the normalised history-dependent equivalent stress \eqref{czc1}, which is a non-linear Abel-type integral operator, implemented before in \cite{MikNam2011AAM} as a (global) material strength condition.
The viscoelasticity is described by a linear Volterra integral operator in time.
The crack starts propagating, breaking the cohesive zone, when the crack tip opening reaches a prescribed critical value. A numerical algorithm for computing the evolution of the crack and cohesive zone in time is discussed along with some numerical results.

As was shown in the paper, the CZ model is applicable only if material parameters, $b$ and $\beta$, of the  history dependent yield condition, based on the power-type durability diagram, are such that $b>0$, $0<\beta<b$. This particularly implies that the CZ model is not applicable for the Robinson-type yield condition, based on the power-type durability diagram.

The CZ model was employed in the plane problem for a single straight crack in an infinite elastic or visco-elastic plane under a homogeneous traction, normal to the crack direction, applied to the plane at infinity at time zero and kept constant in time thereafter.

The numerical results have shown that for both, elastic and visco-elastic materials, there exists a fracture delay time $t_d$, since a remote constant load is applied, during which the cohesive zone grows while the crack does not.

For the growing crack stage, $t>t_d$,  the crack growth rate increases, while the CZ length decreases, with time.
It appeared that in the elastic case,  for some material parameters, there is an unstable crack growth at the onset of crack propagation, followed by stable crack growth. It also causes a jump decrease in the CZ length followed by  a continuous CZ length evolution. However, for other material parameters, no crack instability was detected for the elastic case, implying the stable crack propagation. At the visco-elastic case, the crack propagation was stable for all considered parameters.

The time, when the CZ length decreases to zero seems to coincide with the time when the crack length becomes infinite and can be associated with the complete rupture of the body.
The  rupture time for the visco-elastic case is slightly smaller than that for the purely elastic case.

Implementing different mesh sizes we observed that the solution, normally, converges with the mesh refinement, and we analysed the convergence rates.  An exception is the very slow (if at all) convergence of the CZ tip stress, for some material parameters, which may be a manifestation of a CZ tip stress singularity, at some range of parameters $b$ and $\beta$. Although the square root singularity has been eliminated in the model by the requirement that the corresponding stress intensity factor at the CZ tip is zero, a singularity of a different order can be still present there, however this needs a careful analysis, which is beyond the scope of this paper.

The results presented in the paper particularly show that the the normalised history-dependent equivalent stress \eqref{czc1} is well suitable not only for better approximation of the experimental creep strength data, see \cite{MikNam2011AAM}, but can also be successfully used for numerical solution of some non-stationary problems for bodies under inhomogeneous variable stresses.

As shown in \cite{Mik-ZAMM2000}, the cohesive zone model approach can be also interpreted as a particular non-local approach, cf. \cite{Mikhailov1995-I, Mikhailov1995-II}. In this sense, the CZ model in history-dependent materials presented in this paper is related to another non-local approach based on the Neuber-Novozhilov type stress averaging ahead of the crack under creep or fatigue loading, \cite{MikhNamLiv02, MNGeneva03Creep}.


\end{document}